\documentclass{amsart}

\title{On the algebraic $K$-theory of Witt vectors of finite length}
\author{Vigleik Angeltveit}
\address{Mathematical Sciences Institute \\
Australian National University \\
Australia}
\email{vigleik.angeltveit@anu.edu.au}

\usepackage{amsxtra}
\usepackage{amsfonts}
\usepackage[latin1]{inputenc}
\usepackage{graphicx}
\usepackage{amsmath,amssymb,latexsym,amsthm,mathrsfs}
\usepackage[all]{xy}
\usepackage{pstricks}

\newtheorem{theorem}{Theorem}[section]
\newtheorem{thm}[theorem]{Theorem}
\newtheorem{lemma}[theorem]{Lemma}

\newtheorem{cor}[theorem]{Corollary}
\newtheorem{defn}[theorem]{Definition}

\newtheorem{prop}[theorem]{Proposition}

\newtheorem{remark}[theorem]{Remark}
\newtheorem{example}[theorem]{Example}

\newtheorem{question}[theorem]{Question}

\makeatletter
\let\c@equation\c@theorem
\makeatother
\numberwithin{equation}{section}

\newtheorem{lettertheorem}{Theorem}

  \newcommand{\bC}{\mathbb{C}}   \newcommand{\bF}{\mathbb{F}}          
         \newcommand{\bZ}{\mathbb{Z}}

\newcommand{\sma}{\wedge} %smash product
\newcommand{\holim}{\textnormal{holim}}
\newcommand{\hocolim}{\textnormal{hocolim}}

\newcommand{\coker}{\textnormal{coker}}

\newcommand{\wh}{\widehat}
\newcommand{\wt}{\widetilde}
\newcommand{\xto}{\xrightarrow}

\newcommand{\TR}{\textnormal{TR}}
\newcommand{\TC}{\textnormal{TC}}
\newcommand{\TF}{\textnormal{TF}}

\begin{document}

\begin{abstract}
Let $k$ be a perfect field of characteristic $p$ and let $W_n(k)$ denote the $p$-typical Witt vectors of length $n$. For example, $W_n(\bF_p)=\bZ/p^n$. We study the algebraic $K$-theory of $W_n(k)$, and prove that $K(W_n(k))$ satisfies ``Galois descent''. We also compute the $K$-groups through a range of degrees, and show that the first $p$-torsion element in the stable homotopy groups of spheres is detected in $K_{2p-3}(W_n(k))$ for all $n \geq 2$.
 
\end{abstract}

\maketitle

\section{Introduction}
Let $k$ be a perfect field of characteristic $p$. Then the algebraic $K$-theory of $k$ is well understood, at least after $p$-completion. Indeed, the $p$-completed $K$-theory of $k$ is concentrated in degree $0$.

The situation is far more complicated, but still well understood, if we lift to characteristic $0$ using the Witt vector construction. B\"okstedt and Madsen computed the $p$-completed algebraic $K$-theory of the $p$-adic integers $\bZ_p = W(\bF_p)$ for $p$ odd in \cite{BoMa94}, and of $W(\bF_{p^s})$ (still for $p$ odd) in \cite{BoMa95}. Rognes \cite{Ro99a, Ro99b, Ro99c} computed the $K$-theory of the $2$-adic integers, and Mitchell \cite{Mi02} calculated the $K$-theory of $W(\bF_{2^s})$. Later Hesselholt and Madsen \cite{HeMa03} computed the $K$-theory of complete discrete valuation fields with residue field $k$, at least in odd characteristic.

It follows from the work of Hesselholt and Madsen (\emph{loc.~ cit.}) that $K(W(k))$ satisfies Galois descent. By this we mean that if $k \to k'$ is a $G$-Galois extension of perfect fields of characteristic $p$ with $G$ finite then the canonical map $K(W(k)) \to K(W(k'))^{hG}$ to the homotopy fixed points of $K(W(k'))$ is an equivalence on connective covers after $p$-completion. (This is one version of the Lichtenbaum-Quillen conjecture.) But for $W_n(k)$ for $n<\infty$ the usual tools from algebraic geometry do not work, because $W_n(k)$ is not a regular ring.

Despite considerable effort, very little is known about $K(W_n(k))$. Our first main theorem establishes that $K(W_n(k))$ satisfies Galois descent.

\begin{lettertheorem} \label{t:mainGalois}
Suppose $k \to k'$ is a $G$-Galois extension of perfect fields of characteristic $p$ for a finite group $G$. Then the canonical map
\[
 K(W_n(k)) \to K(W_n(k'))^{hG}
\]
induces an equivalence on connective covers after $p$-completion for any $n$.
\end{lettertheorem}

While the task of understanding $K(W_n(k))$ completely appears insurmountable with current technology, we do have some partial results. The first of those is the following.

\begin{lettertheorem} \label{t:mainorder}
Suppose $k=\bF_q$ is a finite field of characteristic $p$. Then $K_i(W_n(k), (p))$ is a finite $p$-group for all $i \geq 0$ and
\[
 \frac{|K_{2i-1}(W_n(k), (p))|}{|K_{2i-2}(W_n(k), (p))|} = q^{(n-1)i}
\]
for all $i \geq 1$.
\end{lettertheorem}

Combining this with Quillen's calculation of $K(\bF_q)$ we get a similar result for non-relative $K$-theory, see Corollary \ref{c:ordernonrel}.

We can be more explicit in low degrees, explicitly determining the groups up to degree $2p-2$.

\begin{lettertheorem} \label{t:mainlowdeg}
Suppose $k$ is a perfect field of characteristic $p$. Then $K_i(W_n(k), (p))$ for $i \leq 2p-2$ can be described as follows. In odd degree we have
\[
 K_{2i-1}(W_n(k), (p)) \cong W_{(n-1)i}(k)
\]
for $2i-1 \leq 2p-5$ while $K_{2p-3}(W_n(k), (p))$ is the direct sum of $\bZ/p$ and a maximally nontrivial extension of $\coker(\bZ/p \to k)$ by $W_{(n-1)(p-1)-1}(k)$. In even degree we have
\[
 K_{2i}(W_n(k), (p)) \cong \begin{cases} 0 & \textnormal{for $0 \leq 2i \leq 2p-4$} \\
 coker(\phi-1) & \textnormal{for $2i=2p-2$} \end{cases}
\]

Moreover, the unit map from the sphere spectrum sends the first $p$-torsion element $\alpha_1 \in \pi_{2p-3}(S)$ to a generator of $\bZ/p \subset K_{2p-3}(W_n(k), (p))$.
\end{lettertheorem}

\begin{remark}
The above isomorphisms are isomorphisms of abelian groups only. In fact, through the above range of degrees the non-unital graded ring $K_*(W_n(k), (p))$ has trivial multiplication for degree reasons.
\end{remark}

\subsection{Main proof ideas}
This paper uses spectral sequences extensively, both for abstract arguments and for concrete calculations. We strongly encourage the reader to do a few sample calculations on her own to get a feel for how these spectral sequences behave.

We use the cyclotomic trace map \cite{BoHsMa93}
\[
 trc : K(A) \to \TC(A),
\]
which for $A=W_n(k)$ is an equivalence on connective covers after $p$-completion \cite{HeMa97}.

The starting point of our calculation is the topological Hochschild homology of $k$, which looks like $THH(\bF_p)$ ``tensored up'' to $k$. We can bootstrap from that to $THH_*(W_n(k))$ by filtering $W_n(k)$ by powers of $p$. This induces a filtration of $THH(W_n(k))$, and we get a spectral sequence
\[
 E_1^{*,*} = THH_*(k[x]/(x^n)) \Longrightarrow THH_*(W_n(k)).
\]
The $E_1$ term is known by work of Hesselholt and Madsen \cite{HeMa97b}. From this we recover Brun's calculation of $THH_*(\bZ/p^n)$ \cite{Br00}.

This filtration of $THH(W_n(k))$ is $S^1$-equivariant, and as a result we get a corresponding spectral sequence converging to $\TF_*(W_n(k))$. If we use the relative version
\begin{equation} \label{eq:collapsingSSforWnk}
 E_1^{*,*} = \TF_*(k[x]/(x^n), (x)) \Longrightarrow \TF_*(W_n(k), (p)),
\end{equation}
the $E_1$ term is concentrated in odd total degree so the spectral sequence collapses (Corollary \ref{c:collapses}).

The restriction map $R$ does not respect the filtration, so we cannot hope to get such a spectral sequence converging to topological cyclic homology. Instead, the restriction map divides the filtration by $p$, meaning that we have a map
\[
 R : F^s \TF(W_n(k)) \to F^{\lceil s/p \rceil} \TF(W_n(k)),
\]
where $\lceil s/p \rceil$ denotes $s/p$ rounded up to the nearest integer. Expanding on ideas of Brun \cite{Br01} we construct a spectral sequence converging to $\TC_*(W_n(k))$ built from two copies of $\TF_*(W_n(k))$.

The proof of Theorem \ref{t:mainGalois} goes as follows. We first show that $\TF(W_n(k), (p))$ satisfies Galois descent. This follows because $\TF(k[x]/(x^n), (p))$ satisfies Galois descent, plus the collapsing spectral sequence in Equation \ref{eq:collapsingSSforWnk}. Then, because homotopy fixed points commute with homotopy equalizers, the same is true for $\TC$, and the statement for $K$-theory follows by taking connective covers.

The proof of Theorem \ref{t:mainorder} uses the spectral sequence for $\TC$ discussed above. The necessary input is Hesselholt and Madsen's computation of $\TC_*(k[x]/(x^n), (x))$ \cite{HeMa97b}, which implies that for $k=\bF_q$ we have $|\TC_{2i-1}(k[x]/(x^n), (x))|=q^{(n-1)i}$ and $|\TC_{2i-2}(k[x]/(x^n), (x))|=0$.

For Theorem \ref{t:mainlowdeg} we use another idea of Brun \cite{Br01} of comparing with the homotopy orbit spectrum of $THH(W_n(k))$, which is more readily understood. Even though we do not have a map, it is possible to compare $\TC(A)$ to $\Sigma THH(A)_{hS^1}$ because both are the homotopy fiber of a map $\TF(A) \to \TF(A)$. In fact, if $A$ is a complete filtered ring with $I=F^1 A$ we have a spectral sequence converging to $\pi_* \Sigma THH(A,I)_{hS^1}$ with the same $E_1$ term as the above spectral sequence converging to $\TC_*(A,I)$. For $(A,I)=(W_n(k), (p))$ we can compare differentials and extensions through a range of degrees.

\subsection{Relations to previous work}
The calculation of $K_1$ and $K_2$ is classical, and the observation that $K_1(W_n(k)) \cong W_n(k)^{\times}$ behaves differently in characteristic $2$ is of course even more classical (compare $(\bZ/2^n)^\times$ to $(\bZ/p^n)^\times$). Theorem \ref{t:mainlowdeg} can be thought of as an extension of that phenomenon to odd degrees. In characteristic $3$ this was also observed by Geisser \cite{Ge97}, who computed $K_3(W_2(\bF_{3^s}))$ when $(3,s)=1$ and found the extra $\bZ/3$ summand coming from the $3$-torsion in $\pi_3(S)$.

Evens and Friedlander \cite{EvFr82} computed $K_3$ and $K_4$ of $\bZ/p^2$ for $p \geq 5$, but the most general calculation to date, and the only one we know of that goes beyond degree $4$, is due to Brun \cite{Br01} who computed $K_i(\bZ/p^n)$ for $i \leq p-3$.

\subsection{Conventions}
Throughout the paper $k$ will be a perfect field of characteristic $p$ for a fixed prime $p > 0$. We will use the notation $x \doteq y$ to mean $x = \lambda y$ for a unit $\lambda$. We will use $P(x)$, $E(x)$ and $\Gamma(x)$ for a polynomial algebra, exterior algebra, and divided powers algebra, respectively. In $\Gamma(x)$ we will sometimes write $x$ for $\gamma_1(x)$.

\subsection{Acknowledgements}
This paper would never have been started without Mike Hill, with whom I had extensive discussions about the topological Hochschild homology spectral sequence coming from a filtration of a ring. At the time we did not know that Morten Brun had already constructed such a spectral sequence, and we reproved some of his results and did several sample computations together.

In addition I would like to thank Tyler Lawson, Teena Gerhardt, and Lars Hesselholt for helpful conversations. I would also like to thank an anonymous referee for constructive criticism and for spotting a serious mistake in an earlier version of this paper.

Because of its long gestation period this work has been supported by several grants: An NSF All-Institutes Postdoctoral Fellowship administered by the Mathematical Sciences Research Institute through its core grant DMS-0441170, NSF grant DMS-0805917, and Australian Research Council Discovery Grant No.\ DP120101399.

\section{A topological Hochschild homology spectral sequence}
In this section we study a spectral sequence
\[
 E_1^{s,t} = \pi_{s+t} THH(Gr A; s) \Longrightarrow \pi_{s+t} THH(A)
\]
associated to a filtration of a ring $A$. We will call this the ``filtered ring spectral sequence'' for $THH$. The existence of this spectral sequence was first noted by Brun \cite{Br00}, though he only used it in an indirect way in his computation of $THH_*(\bZ/p^n)$. We will demonstrate that this spectral sequence is an excellent tool for computations by simplifying and extending known calculations of $THH$.

For conventions and standard results about spectral sequences, see \cite{Bo99}. Most of the spectral sequences in this paper will be conditionally convergent. If the spectral sequence satisfies some Mittag-Leffler condition it converges strongly. This is typically easy to verify, in most of our examples it follows because the $E_1$-term is finite (or has finite length) over $k$ in each bidegree. Because of the large number of spectral sequences appearing we will not discuss convergence in each case.

\subsection{A Hochschild homology spectral sequence} \label{s:HH}
We start with Hochschild homology, which is easier, in order to introduce some key ideas. Recall that for a ring $A$, the Hochschild homology $HH_*(A)$ is the homology of a chain complex $C_*(A)$ with $A^{\otimes q+1}$ in degree $q$ and
\begin{multline*}
d(a_0 \otimes \ldots \otimes a_q) = \sum_{0 \leq i \leq q-1} (-1)^i a_0 \otimes \ldots \otimes  a_i a_{i+1} \otimes \ldots \otimes a_q \\
+ (-1)^q a_q a_0 \otimes a_1 \otimes \ldots \otimes a_{q-1}.
\end{multline*}
It can also be described as the homology of the cyclic bar construction $B^{cy}_\otimes(A)$.

If $A$ is a DGA rather than a ring the above definition yields a bicomplex, and the Hochschild homology of $A$ is the homology of the associated total complex. As always we follow the usual sign rule, multiplying by $(-1)$ whenever we move two things (elements, or operators like $d$) of odd degree past each other.

If the ring $A$ is not flat as a $\bZ$-module, the above definition of Hochschild homology will give the ``wrong'' result because we were using underived tensor products while we should have been using derived tensor products. Shukla homology is the appropriate version of Hochschild homology with derived tensor products everywhere. Alternatively, we can replace $A$ by a DGA $\widetilde{A}$ which is degree-wise flat and satisfies $H_*(\widetilde{A}) \cong A$. One way to do this is to put a model category structure on the category of DGAs \cite{Ja97} and let $\wt{A}$ be a cofibrant replacement of $A$. Then $\wt{A}$ is degree-wise projective, hence flat. It is standard that $HH_*(\widetilde{A})$ is independent of the replacement $\widetilde{A}$ and agrees with the Shukla homology of $A$.

An alternative description of $HH_*(A)$, at least when $A$ is flat, is as the homology of the derived tensor product $A \otimes_{A \otimes A^{op}}^L A$, or as $Tor^{A \otimes A^{op}}_*(A, A)$. The equivalence between these last two definitions follows by replacing one of the $A$'s by the $2$-sided bar construction $B(A,A,A)$, which is a cofibrant replacement of $A$ as an $A$-bimodule.

\begin{example} \label{ex:HHofZ/p}
Let $A=\bZ/p$. Because $\bZ/p$ is not projective, we replace it by the DGA $\widetilde{\bZ/p}$ defined as follows. As a chain complex it is given by
\[
 \bZ\{a\} \overset{d}{\to} \bZ\{1\}.
\]
Here the generator $a$ is in degree $1$, the generator $1$ is in degree $0$, and the boundary map is given by $d(a)=p \cdot 1$. There is only one way to define the multiplication: we must have $a^2=0$.

We can then compute the Hochschild homology of $\widetilde{\bZ/p}$. We find that $\widetilde{\bZ/p} \otimes \widetilde{\bZ/p} \cong \widetilde{\bZ/p} \otimes E_\bZ(b)$, where $E_\bZ(b)$ denotes an exterior algebra (over $\bZ$) on a generator $b$ in degree $1$. We can take $b=1 \otimes a - a \otimes 1$. It follows that $H_*(\widetilde{\bZ/p} \otimes \widetilde{\bZ/p}) \cong E_{\bZ/p}(b)$, an exterior algebra over $\bZ/p$ on one generator $b$ in degree $1$. It then follows from standard homological algebra that
\[
 HH_*(\widetilde{\bZ/p}) \cong Tor_*^{E_{\bZ/p}(b)}(\bZ/p, \bZ/p) \cong \Gamma_{\bZ/p}(\mu_0)
\]
is a divided powers algebra over $\bZ/p$ on a class $\mu_0$ in degree $2$. The class $\mu_0$ is represented by $1 \otimes a - a \otimes 1$ in the Hochschild complex; it has degree $2$ because it has internal degree $1$ and Hochschild degree $1$.
\end{example}

Suppose $A = \bigoplus_{i \in \bZ} A_i$ is a graded ring. In the examples this grading will usually be independent of the homological grading. Then we get a splitting of the Hochschild homology of $A$.

\begin{lemma} \label{l:HHsplitting}
Suppose $A$ is a graded ring. Then the Hochschild homology $HH_*(A)$ of $A$ splits as a direct sum
\[
 HH_*(A) \cong \bigoplus_s HH_*(A;s),
\]
where $HH_*(A;s)$ is the homology of the subcomplex of $C_*(A)$ of internal degree $s$. Here we give $a_0 \otimes \ldots \otimes a_q$ in $C_q(A)$, with each $a_i$ homogeneous, internal degree $|a_0|+\ldots+|a_q|$.
\end{lemma}

\begin{proof}
This is clear, because the Hochschild differential preserves the internal degree.
\end{proof}

Now suppose $A$ is a complete filtered ring. By this we mean that $A$ comes with a decreasing filtration
\[
 \ldots \subset F^{s+1} A \subset F^s A \subset \ldots \subset F^0 A=A.
\]
We assume the filtration is compatible with the multiplicative structure, meaning that the multiplication on $A$ induces maps $F^i A \otimes F^j A \to F^{i+j} A$. Complete means that the canonical map $A \to \lim_s A/F^s A$ is an isomorphism. The canonical example comes from an ideal $I \subset A$. If $A$ is $I$-complete then $F^s A=I^s A$ defines a complete filtration on $A$. Let $Gr^i A=F^i A/F^{i+1} A$ and let $Gr A = \bigoplus_i Gr^i A$. Then $Gr A$ is a graded ring, and we can compute $HH_*(Gr A)$ as above.

To get homotopically meaningful results we do the following. First, we replace $A$ by a projective DGA $\wt{A}$ and define $F^s \wt{A}$ by basechange along $F^s A \to A$. Then each $F^s \wt{A}$ is degreewise projective, hence flat, and the multiplication map on $\wt{A}$ induces maps $F^i \wt{A} \otimes F^j \wt{A} \to F^{i+j} \wt{A}$. ($F^i \wt{A} \otimes F^j \wt{A}$ maps to both $\wt{A}$ and to $F^{i+j} A$, hence maps to the pullback $F^{i+j} \wt{A}$.) Second, we define
\[
 \wt{F}^s \wt{A} = \underset{t \geq s}{\hocolim} F^t \wt{A},
\]
which is chain homotopic to $F^s \wt{A}$. The multiplication map $F^i \wt{A} \otimes F^j \wt{A} \to F^{i+j} \wt{A}$ then induces a map $\wt{F}^i \wt{A} \otimes \wt{F}^j \wt{A} \to \wt{F}^{i+1} \wt{A}$, and we define
\[
 \wt{Gr}^s \wt{A} = \wt{F}^s \wt{A}/\wt{F}^{s+1} \wt{A}.
\]
Then $\wt{Gr}^s \wt{A}$ is chain homotopic to $F^s \wt{A}/F^{s+1} \wt{A}$, but has the advantage that each $\wt{Gr}^s \wt{A}$ is flat. The direct sum
\[
 \wt{Gr} \wt{A} = \bigoplus_{s \geq 0} \wt{Gr}^s \wt{A}
\]
is a graded ring, so we get a splitting of $HH_*(\wt{Gr} \wt{A})$ as in Lemma \ref{l:HHsplitting}.

Next we define a corresponding filtration of $C_*(A)$, or rather of $C_*(\wt{A})$. To get a homotopically meaningful result we do this by defining
\[
 \wt{F}^s C_q(\wt{A}) = \underset{s_0+\ldots+s_q \geq s}{\hocolim} F^{s_0} \wt{A} \otimes \ldots \otimes F^{s_q} \wt{A}.
\]
Then $\wt{F}^0 C_q(\wt{A})$ is chain homotopy equivalent to $C_q(\wt{A})$, so we can think of
\[
 \ldots \subset \wt{F}^{s+1} C_q(\wt{A}) \subset \wt{F}^s C_q(\wt{A}) \subset \ldots \subset \wt{F}^0 C_q(\wt{A})
\]
as a filtration of $C_q(\wt{A})$. Then we define
\[
 \wt{Gr}^s C_q(\wt{A}) = \wt{F}^s C_q(\wt{A})/\wt{F}^{s+1} C_q(\wt{A})
\]
and
\[
 \wt{Gr} C_q(\wt{A}) = \bigoplus_{s \geq 0} \wt{Gr}^s C_q(\wt{A}).
\]

The Hochschild differential maps the diagram defining $\wt{F}^s C_q(\wt{A})$ to the diagram defining $\wt{F}^s C_{q-1}(\wt{A})$, so it induces a map $\wt{F}^s C_q(\wt{A}) \to \wt{F}^s C_{q-1}(\wt{A})$. Hence we have a filtration of $\wt{F}^0 C_*(\wt{A})$, and this filtration gives us a spectral sequence. We can identify the $E_1$ term of this spectral sequence as follows:

\begin{thm}
Suppose $A$ is a complete filtered ring with associated graded $Gr A$. Then there is a weakly convergent spectral sequence
\[
 E_1^{s,t} = HH_{s+t}(\widetilde{Gr A};s) \Longrightarrow HH_{s+t}(\widetilde{A}).
\]
The differential $d_r$ has bidegree $(r,-r-1)$. If $A$ is commutative this is an algebra spectral sequence.
\end{thm}

We only get weak convergence in general because we have no guarantee that $\wt{F}^\infty C_q(\wt{A})$ is contractible even if $F^\infty A=0$.

\begin{proof}
It follows that with the above filtration of $\wt{F}^0 C_q(\wt{A})$, $\wt{Gr}^s C_q(\wt{A})$ is chain homotopy equivalent to
\[
 \bigoplus_{s_0+\ldots+s_q=s} \wt{Gr}^{s_0} \wt{A} \otimes \ldots \otimes \wt{Gr}^{s_q}(\wt{A}).
\]
This computes the $E_1$ term of the spectral sequence.

If $A$ is commutative then we can choose $\wt{A}$ to be commutative as well, and each $\wt{Gr} C_q(\wt{A})$ becomes a commutative DGA. Then the shuffle product induces a multiplication on the spectral sequence which makes it an algebra spectral sequence in the standard way.
\end{proof}

Note that being an algebra spectral sequence means that the differentials are graded derivations: $d_r(xy)=d_r(x) y + (-1)^{|x|}x d_r(y)$ for all $r \geq 1$ and $x,y \in E_r^{*,*}$.

Next we look at some examples to show that this spectral sequence can be used quite effectively. We fix a perfect field $k$ of characteristic $p$. Then by a generalization of Example \ref{ex:HHofZ/p} (use \cite[Lemma 5.5]{HeMa97}) the Hochschild homology of $\widetilde{k}$ is a divided powers algebra over $k$ on one generator $\mu_0$ in degree $2$.

\begin{example} \label{ex:HHofWk}
First we consider $W(k)$ filtered by powers of $p$. Then the associated graded is $Gr W(k) \cong k[x]$, and we have a spectral sequence
\[
 E_1^{s,t} = HH_{s+t}(\widetilde{k}[x]; s) \Longrightarrow HH_*(\widetilde{W(k)}).
\]
A priori this spectral sequence only converges weakly because $\wt{F}^\infty C_q(\wt{A})$ is a large rational vector space. But the $p$-completion of $\wt{F}^\infty C_q(\wt{A})$ is trivial, so the spectral sequence converges strongly to the $p$-completion of $HH_*(\wt{W(k)})$.

We find that
\[
 E_1^{*,*} = HH_*(\widetilde{k}[x]) \cong P(x) \otimes E(\sigma x) \otimes \Gamma(\mu_0),
\]
where $\mu_0$ comes from $HH_*(\widetilde{k})$ and $\sigma x$ is represented by $1 \otimes x - x \otimes 1$. Here the tensor product is over $k$. This is bigraded, with $|\mu_0|=(0,2)$, $|x|=(1,-1)$ and $|\sigma x|=(1,0)$.

We have an immediate differential
\[
 d_1(\gamma_j(\mu_0)) = \gamma_{j-1}(\mu_0) \sigma x
\]
for each $j \geq 1$, leaving
\[
 E_2^{*,*} = E_\infty^{*,*} = P(x)
\]
concentrated in homological degree $0$.

If in addition we use that there is a comultiplication on $E_1^{*,*}$ with $\psi(\gamma_j(\mu_0)) = \sum_{a+b=j} \gamma_a(\mu_0) \otimes \gamma_b(\mu_0)$ as in \cite{AnRo05} we find that the $d_1$-differential is generated by the single differential $d_1(\mu_0) = \sigma x$. Since $x$ represents multiplication by $p$, this recovers (at least up to $p$-completion) the classical result that $HH_0(\widetilde{W(k)})=W(k)$ and $HH_i(\widetilde{W(k)})=0$ for $i > 0$.
\end{example}

\begin{example} \label{ex:HHofWnk}
Next we consider $W_n(k)$ filtered by powers of $p$. Then the associated graded is $Gr W_n(k) = k[x]/(x^n)$, and in this case there are no convergence problems. Let
\[
 E_0^{*,*} = P_n(x) \otimes E(\sigma x) \otimes \Gamma(x_n) \otimes \Gamma(\mu_0),
\]
where the new generator $x_n$ has bidegree $|x_n|=(n,2-n)$. Now define a differential $d_0$ on $E_0^{*,*}$, generated multiplicatively by $d_0(\gamma_j(x_n))=nx^{n-1} \gamma_{j-1}(x_n) \sigma x$ for $j \geq 1$. Then
\[
 E_1^{*,*} = HH_*(\widetilde{k}[x]/(x^n)) \cong H_*(E_0^{*,*}, d_0).
\]
If $p$ divides $n$ then $d_0=0$, and $E_1^{*,*}=E_0^{*,*}$ with a $d_1$-differential generated multiplicatively by $d_1(\gamma_j(\mu_0))=\gamma_{j-1}(\mu_0) \sigma x$ for $j \geq 1$, leaving
\[
 E_2^{*,*} = E_\infty^{*,*} = P_n(x) \otimes \Gamma(x_n).
\]
This is the associated graded of $HH_*(\widetilde{W_n(k)}) \cong W_n(k) \otimes \Gamma(x_n)$. As above, if we use that there is a comultiplication on $E_1^{*,*}$ with $\psi(\gamma_j(\mu_0)) = \sum_{a+b=j} \gamma_a(\mu_0) \otimes \gamma_b(\mu_0)$ we can say that the $d_1$-differential is generated by the single differential $d_1(\mu_0) = \sigma x$.

If $p$ does not divide $n$ then the $E_1$-term is somewhat smaller. We still have a $d_1$-differential generated by $d_1(\mu_0)=\sigma x$, but now the $E_2$-term is somewhat larger. In this case we also have $d_2$-differentials
\[
 d_2(\gamma_i(x_n) x^{n-1} \gamma_j(\mu_0)) \doteq \gamma_{i+1}(x_n) \gamma_{j-2}(\mu_0) \sigma x
\]
for $j \geq 2$. We illustrate this with the zig-zag
\[ \xymatrix@R-18pt{
 \gamma_i(x_n) x^{n-1} \gamma_j(\mu_0) \ar@{|->}[dr] & \\
 & \gamma_i(x_n) x^{n-1} \gamma_{j-1}(\mu_0) \sigma x \\
 \gamma_{i+1}(x_n) \gamma_{j-1}(\mu_0) \ar@{|->}[ur] \ar@{|->}[dr] & \\
 & \gamma_{i+1}(x_n) \gamma_{j-2}(\mu_0) \sigma x
} \]
This leaves
\[
 E_3^{*,*} = E_\infty^{*,*} = P_n(x)\{1\} \oplus \bigoplus_{j \geq 1} \big( k\{x^{n-1} \mu_0 \gamma_{j-1}(x_n)\} \oplus P_{n-1}(x)\{x \gamma_j(x_n)\} \big).
\]
There is a hidden multiplication by $p$ extension, so again we recover that
\[
 HH_*(\widetilde{W_n(k)}) \cong W_n(k) \otimes \Gamma(x'_n).
\]
Now $\gamma_j(x'_n)$ is represented by $x^{n-1} \mu_0 \gamma_{j-1}(x_n)$, while $p\gamma_j(x'_n)$ is represented by $x \gamma_j(x_n)$.
\end{example}

\begin{remark} \label{r:filter_away}
Note that in the above example the case $p \nmid n$ is more complicated. It is possible to filter away this added complexity, as follows. In the Hochschild chain complex $C_*(\widetilde{W_n(k)})$, introduce a third grading by giving the class representing $\gamma_j(x_n)$ degree $-j$ with associated graded $\widehat{Gr}C_*(\widetilde{W_n(k)})$. Then we get a spectral sequence $(\wh{Gr} E_r^{*,*}, \wh{d}_r)$ converging to $\widehat{Gr} HH_*(\widehat{W_n(k)})$. The associated graded $\widehat{Gr} E_0^{*,*}$ is the ring $E_0^{*,*}$ above, now trigraded. Then we get the same $d_1$-differential as in the case $p \mid n$, at which point the spectral sequence once again collapses. We now have another spectral sequence
\[
 E_1^{*,*} = W_n(k) \otimes \Gamma(x_n) \Longrightarrow HH_*(\widetilde{W_n(k)}),
\]
which collapses at the $E_1$-term, giving us the desired result without having to compute higher differentials.
\end{remark}

\subsection{Topological Hochschild homology} \label{s:THH}
For a naive definition of $THH$ we have a wide choice of frameworks with which to work. For example, we could define $THH(A)$ as the geometric realization of a simplicial spectrum with $q \mapsto A^{\sma q+1}$, the $(q+1)$-fold smash product of $A$ with itself. But to build $THH(A)$ as a cyclotomic spectrum (see Section \ref{ss:fp} below for the definition of a cyclotomic spectrum) we need a more sophisticated definition. A variant of this definition goes back to B\"okstedt \cite{Bo1}, see also \cite{HeMa97}.

Since this technology is well established, we will be brief. Our definition will take as input a symmetric ring spectrum and give us back an orthogonal $S^1$-spectrum; our definition should be compared to that given by Hesselholt and Madsen in \cite{HeMa97}. See \cite{ABGHLM} for a modern definition of $THH$ which takes place entirely in orthogonal spectra.

Let $A$ be a symmetric ring spectrum in the sense of \cite{HSS}, but with topological spaces instead of simplicial sets. For convenince we will assume that all our symmetric ring spectra are strictly connective and convergent, i.e., $A_n$ is $(n-1)$-connected and there exists a sequence of nondecreasing integers $\alpha_n$ tending to infinity such that $\Sigma A_n \to A_{n+1}$ is $(n+\alpha_n)$-connected. 

If $A$ is a ring, we can regard $A$ as a symmetric ring spectrum by setting $A_i=K(A,i)$ for a particular choice of $K(A,i)$, see \cite[Example 1.2.5]{HSS}. For each simplicial degree $q$ and finite-dimensional real inner product space $V$ we can consider the space
\[
 THH(A)_q(V) = \underset{I^{q+1}}{\hocolim} \Omega^{i_0+\ldots+i_q} (A_{i_0} \sma \ldots \sma A_{i_q} \sma S^V).
\]
Here $I$ is the category whose objects are ${\bf n} = \{1,\ldots,n\}$ for $n \geq 0$ and whose morphisms are all injective maps. By varying $n$ we get an orthogonal spectrum (in the sense of \cite{MMSS}) $THH(A)_q$ for each $q$, and by varying $q$ we get a simplicial orthogonal spectrum. We define the orthogonal spectrum $THH(A)$ as the geometric realization of this simplicial orthogonal spectrum. Because it is the geometric realization of a cyclic object we get an $S^1$-action.

It is possible to define an orthogonal $S^1$-spectrum as an orthogonal spectrum with an $S^1$-action. To distinguish between ``naive'' and ``genuine'' equivariant spectra one can vary the model structure. If we wish to evaluate $THH(A)$ on an $S^1$-representation $V$, we can define $THH(A)(V)$ as above with $S^1$ acting diagonally.

We will discuss the model structure on orthogonal $S^1$-spectra in Section \ref{s:tracemethod} below.

In unpublished work \cite{Bo2}, B\"okstedt computed $THH(\bF_p)$ and $THH(\bZ)$, and in \cite{HeMa97} Hesselholt and Madsen extended the first of these calculations to $THH(k)$ for any perfect field $k$ of characteristic $p$. They found that
\[
 \pi_* THH(k) \cong P(\mu_0),
\]
a polynomial algebra over $k$ on one variable $\mu_0$ in degree $2$. Here $\mu_0$ is represented by the image of $\bar{\tau}_0 \in \pi_1(H\bF_p \sma H\bF_p)$ in $\pi_1(Hk \sma Hk)$, where $\tau_0$ is the mod $p$ Bockstein and $\bar{\tau}_0=-\tau_0$ is its conjugate. If we put a CW structure on $H\bF_p$ with one $0$-cell denoted $1$ and one $1$-cell denoted $a$ then $\mu_0$ is represented by $1 \sma a - a \sma 1$ in the cellular chain complex of $H\bF_p \sma H\bF_p$. It follows that the class $\mu_0$ maps to the class with the same name in $HH_2(\wt{k})$.

Let $\nu_p(i)$ denote the $p$-adic valuation of $i$. We will see in Example \ref{ex:THHofWk} below that
\[
 \pi_j THH(W(k)) \cong \begin{cases} W(k) & \textnormal{if $j=0$} \\ W_{\nu_p(i)}(k) & \textnormal{if $j=2i-1$ is odd} \\ 0 & \textnormal{if $j \neq 0$ is even} \end{cases}
\]
These spectra are sometimes easier to understand if we use mod $p$ coefficients. Let $V(0)$ denote the mod $p$ Moore spectrum. Then
\[
 V(0)_* THH(W(k)) \cong E(\lambda_1) \otimes P(\mu_1),
\]
where the ground ring is $k$ and $|\lambda_1|=2p-1$, $|\mu_1|=2p$.

We can then recover $THH_*(W(k))$ by running the Bockstein spectral sequence
\[
 E_1^{*,*} = V(0)_* THH(W(k))[v_0] \Longrightarrow THH_*(W(k)).
\]
This spectral sequence is generated multiplicatively by the differentials
\[
 d_{j+1}(\mu_1^{p^j}) = v_0^{j+1} \mu_1^{p^j-1} \lambda_1
\]
for $j \geq 0$. If in addition we use the ``Leibniz rule'' $d_{j+1}(y^p)=v_0 y^{p-1} d_j(y)$ then the Bockstein spectral sequence is generated by the single differential $d_1(\mu_1)=v_0 \lambda_1$.

\begin{remark} \label{r:Leibniz}
The ``Leibniz rule'' in the Bockstein spectral sequence going from mod $p$ homology to integral homology is discussed in \cite[Proposition 6.8]{Ma70}; at $p=2$ there is a correction term for $d_2$ but otherwise it holds. While we have mod $p$ and integral homotopy instead of homology, a similar result holds. The correction term for $d_2$ at $p=2$ is $Q^4(\lambda_1$), and an explicit computation shows that this is indeed $0$.
\end{remark}

Returning to the general theory, suppose $A$ is a graded ring. Then we get a splitting of $THH(A)$ into homogeneous pieces in the same way as for Hochschild homology.

\begin{lemma}
Suppose $A$ is a graded ring or symmetric ring spectrum. Then
\[
 THH(A) \cong \bigvee_s THH(A;s),
\]
where $THH(A;s)$ is the geometric realization of the subcomplex $THH(A;s)_\bullet$ of internal degree $s$.
\end{lemma}

\begin{proof}
Write $Gr^s A$ for the $s$'th graded piece of $A$, and define
\[
 THH(A;s)_q(V) = \bigvee_{s_0+\ldots+s_q=s} \underset{I^{q+1}}{\hocolim} \, \Omega^{i_0+\ldots+i_q}(Gr^{s_0} A_{i_0} \sma \ldots \sma Gr^{s_q} A_{i_q} \sma S^V).
\]
The face and degeneracy maps respect this splitting, hence we get a corresponding splitting after geometric realization.
\end{proof}

\subsection{A topological Hochschild homology spectral sequence} \label{ss:THHSS}
Now suppose $A$ is a complete filtered ring or symmetric ring spectrum. Recall that we are assuming that $A$ is strictly connective and convergent. We will assume that each $F^s A$ is strictly connective and convergent as well. We can then define a corresponding filtration on $THH(A)$, by setting
\[
 F^s THH(A)_q(V) = \underset{i_0+\ldots+i_q \geq s}{\hocolim} \underset{I^{q+1}}{\hocolim} \, \Omega^{i_0+\ldots+i_q} \big( F^{s_0} A_{i_0} \sma \ldots \sma F^{s_q} A_{i_q} \sma S^V \big).
\]
We note that there is an induced map $F^s THH(A)_q(V) \sma S^W \to F^s THH(A)_q(V \oplus W)$, so each $F^s THH(A)_q$ is indeed an orthogonal spectrum. Moreover, this filtration is compatible with the face and degeneracy maps, so we can define $F^s THH(A)$ as the geometric realization of the simplicial orthogonal spectrum $q \mapsto F^s THH(A)_q$. We can now prove a result due to Brun, stated in more modern language.

\begin{thm}[Brun \cite{Br00}] \label{t:THHSS}
Suppose $A$ is a complete filtered ring or symmetric ring spectrum with associated graded $Gr A$. Then there is a weakly convergent spectral sequence
\[
 E_1^{s,t} = THH_{s+t}(Gr A;s) \Longrightarrow THH_{s+t}(A).
\]
If $A$ is commutative this is an algebra spectral sequence.
\end{thm}

\begin{proof}
(Compare \cite[Proposition 5.2.2]{Br00}.) We have
\[
 F^s THH(A)_q(V) = \underset{s_0+\ldots+s_q \geq s}{\hocolim} \underset{I^{q+1}}{\hocolim} \Omega^{i_0+\ldots+i_q} \big( F^{s_0} A_{i_0} \sma \ldots \sma F^{s_q} A_{i_q} \sma S^V \big).
\]
We can exchange the two homotopy colimits and consider
\[
 \underset{s_0+\ldots+s_q \geq s}{\hocolim} \Omega^{i_0+\ldots+i_q} \big( F^{s_0} A_{i_0} \sma \ldots \sma F^{s_q} A_{i_q} \sma S^V \big).
\]
Because $A$ is connective, it follows that the map
\begin{multline*}
 \underset{s_0+\ldots+s_q \geq s}{\hocolim} \Omega^{i_0+\ldots+i_q} \big( F^{s_0} A_{i_0} \sma \ldots \sma F^{s_q} A_{i_q} \sma S^V \big) \\
 \to \Omega^{i_0+\ldots+i_q} \big( \underset{s_0+\ldots+s_q \geq s}{\hocolim} F^{s_0} A_{i_0} \sma \ldots \sma F^{s_q} A_{i_q} \sma S^V \big)
\end{multline*}
becomes increasingly connective as the dimension of $V$ goes to infinity.

It follows that the filtration quotient $F^s THH(A)_q(V)/F^{s+1} THH(A)_q(V)$ is equivalent to
\[
 \Omega^{i_0+\ldots+i_q} \big( \bigvee_{s_0+\ldots+s_q=s} Gr^{s_0} A_{i_0} \sma \ldots \sma Gr^{s_q} A_{i_q} \sma S^V \big)
\]
through a range that goes to infinity as the dimension of $V$ goes to infinity.

If $A$ is commutative the shuffle product induces an algebra structure on the spectral sequence in the usual way.
\end{proof}

\begin{remark}
To get a multiplication on the spectral sequence it suffices to assume that $A$ is an $E_2$ ring spectrum. This is related to how $THH(A)$ is an $S$-algebra as long as $A$ is an $E_2$ ring spectrum, see \cite{BrFiVo07}. We omit the details, as we will not need them.
\end{remark}

\subsection{Sample computations} \label{s:ex}
In this section we use Theorem \ref{t:THHSS} to compute $THH(A)$ in some examples.

\begin{example} \label{ex:THHofWk}
We start by computing $THH_*(W(k))$ from $THH_*(k[x])$. We find that
\[
 E_1^{*,*} = THH_*(k[x]) \cong P(x) \otimes E(\sigma x) \otimes P(\mu_0),
\]
where $\mu_0$ comes from $THH_*(k)$ and the tensor product is over $k$. The only difference from the Hochschild homology calculation in Example \ref{ex:HHofWk} is that here $\mu_0$ is a polynomial generator rather than a divided powers generator.

We have an immediate differential
\[
 d_1(\mu_0)=\sigma x,
\]
because $\mu_0$ is represented by $1 \otimes \bar{\tau}_0$ where $\tau_0$ is the mod $p$ Bockstein and $\sigma x$ is represented by $1 \otimes x$. Hence
\[
 E_2^{*,*} = P(x) \otimes E(\lambda_1) \otimes P(\mu_1),
 \]
where $\mu_1=\mu_0^p$ and $\lambda_1 = \mu_0^{p-1} \sigma x$. Next we use the Leibniz rule to get a differential $d_2(\mu_1)=x\lambda_1$, so
\[
 E_3^{*,*} =  P(x) \otimes E(\lambda_2) \otimes P(\mu_2) \oplus \{\textnormal{$x$-torsion}\}.
\]
In general we find that
\[
 E_{r+1}^{*,*} = P(x) \otimes E(\lambda_r) \otimes P(\mu_r) \oplus \{\textnormal{$x$-torsion}\},
\]
where $\mu_r = \mu_{r-1}^p$ and $\lambda_r = \mu_{r-1}^{p-1} \lambda_{r-1}$, and we recover $THH_*(W(k))$. Note that the $E_2$-term of this spectral sequence is isomorphic to the $E_1$-term of the Bockstein spectral sequence which computes $THH_*(W(k))$ from $V(0)_* THH(W(k))[v_0]$, discussed immediately before Remark \ref{r:Leibniz}.
\end{example}

\begin{example} \label{ex:THHofW_nk}
Next we compute $THH_*(W_n(k))$, starting from $THH_*(k[x]/(x^n))$. As for Hochschild homology, the calculation is easier if $p \mid n$. Let
\[
 E_0^{*,*} = P_n(x) \otimes E(\sigma x) \otimes \Gamma(x_n) \otimes P(\mu_0) 
\]
and define a differential $d_0$ on $E_0$ by $d_0(x_n)=n x^{n-1} \sigma x$. Then
\[
 E_1^{*,*} = THH_*(k[x]/(x^n)) \cong H_*(E_0^{*,*},d_0).
\]

First suppose $p \mid n$. Then $d_0=0$, so $E_1^{*,*}=E_0^{*,*}$ and we get the same differentials
\[
 d_{j+1}(\mu_0^{p^j}) = x^j \mu_0^{p^j-1} \sigma x
\]
as for $W(k)$, for $0 \leq j \leq n-1$.

Next suppose $p \nmid n$. Then, just as in the computation of $HH_*(W_n(k))$, this moves the differentials around. The end result is that the $E_\infty$ term is isomorphic to what we get in the case $p \mid n$, except that we have a class $\gamma_{m-1}(x_n) x^{n-1} \mu_0$ instead of $\gamma_m(x_n)$ for each $m \geq 0$.

Alternatively, we can follow the approach in Remark \ref{r:filter_away}. As for Hochschild homology, we introduce another filtration on $THH(W_n(k))$ so that the associated graded is the ring $E_0^{*,*}$ above, now trigraded. This reduces the case $p \nmid n$ to the case $p \mid n$. This proves the following:

\begin{thm} \label{t:THHofW_nk}
The ring $THH_*(W_n(k))$ is the homology of the DGA
\[
 A = \Gamma_{W_n(k)}(x'_n) \otimes E_{W_n(k)}(\lambda_0) \otimes P_{W_n(k)}(\mu_0)
\]
where $\lambda_0$ is represented by $\sigma x$, $x'_n$ is represented by $x_n-np^{n-1} \mu_0$, and the tensor product is over $W_n(k)$. The differential is multiplicatively generated by $d(\mu_0) = \lambda_0$. The homology is given by
\begin{eqnarray*}
 THH_{2i}(W_n(k)) & \cong & \bigoplus_{0 \leq j \leq i} W_{\max(\nu_p(j), n)}(k) \quad \textnormal{for } i \geq 0 \\
 THH_{2i-1}(W_n(k)) & \cong & \bigoplus_{1 \leq j \leq i} W_{\max(\nu_p(j), n)}(k) \quad \textnormal{for } i \geq 1
\end{eqnarray*}
\end{thm}

This recovers Brun's calculation of $THH_*(\bZ/p^n)$ from \cite{Br00}. We note that the first nonzero odd group is $THH_{2p-1}(W_n(k)) \cong k$, and that the canonical map $THH(W(k)) \to THH(W_n(k))$ maps $THH_{2p-1}(W(k)) \cong k$ isomorphically onto this $k$.
\end{example}

We include one more example. This next example will not be used in the rest of the paper.

\begin{example}
Consider the Adams summand $\ell$ of connective $p$-local complex $K$-theory $ku_{(p)}$. We filter this by powers of $v_1$:
\[
 \ldots \to \Sigma^{(n+1)(2p-2)} \ell \to \Sigma^{n(2p-2)} \ell \to \ldots \to \ell.
\]
This filtration is multiplicative, and the associated graded is
\[
 Gr \ell \cong H\bZ_{(p)}[v_1],
\]
where $|v_1|=2p-2$.

Now, consider the resulting spectral sequence with mod $p$ coefficients. We find that
\[
 E_1^{*,*} = V(0)_* THH(\bZ_{(p)}[v_1]) \cong E(\lambda_1) \otimes P(\mu_1) \otimes P(v_1) \otimes E(\sigma v_1),
\]
and there is an immediate differential $d_1(\mu_1) = \sigma v_1$, leaving us with
\[
 E_2^{*,*} = P(\mu_2) \otimes E(\lambda_1,\lambda_2) \otimes P(v_1).
\]
Here $\mu_2=\mu_1^p$ and $\lambda_2 = \mu_1^{p-1} \sigma v_1$. This coincides with the $E_1$-term of the $v_1$-Bockstein spectral sequence considered in \cite{McSt93}.

This spectral sequence is also interesting with integral coefficients. Recall from \cite{AnHiLa10} that in $THH_*(\ell)$ there is an infinite $v_1$-tower on $\lambda_1$ which becomes increasingly $p$-divisible. In $THH_{2p-1}(\bZ_{(p)}[v_1])$ there is a $\bZ/p$ generated by $\lambda_1$ and a $\bZ_{(p)}$ generated by $\sigma v_1$, and there is a nontrivial extension $p \cdot \lambda_1 = \sigma v_1$ in $THH_*(\ell)$. Hence the class $\lambda_1$ is $1/p$ times a naturally defined class.

We have not attempted to understand the general behavior of the spectral sequence $THH_*(\bZ_{(p)}[v_1]) \Longrightarrow THH_*(\ell)$, though it is interesting that with the two spectral sequences in \cite{AnHiLa10} we now have three spectral sequences all converging to $THH_*(\ell)$.
\end{example}

\subsection{Relative topological Hochschild homology}
Given an ideal $I \subset A$, we define $THH(A,I)$ to be the homotopy fiber of the canonical map $THH(A) \to THH(A/I)$. If $A$ is a filtered ring as above, $I=F^1 A$ becomes an ideal and $THH(A,I)$ is homotopy equivalent to $F^1 THH(A)$. Hence we get a spectral sequence converging to $THH_*(A,I)$ with $E_1$ term
\[
 E_1^{s,t} = \begin{cases} THH_{s+t}(Gr A; s) \quad & s \geq 1 \\ 0 \quad & s=0 \end{cases}
\]
We analyze the effect of removing filtration $0$ in some examples.

\begin{example}
Consider $THH(W(k), (p))$ with $W(k)$ filtered by powers of $p$. Then we have a spectral sequence
\[
 E_1^{*,*} = ker \big( P(x) \otimes E(\sigma x) \otimes P(\mu_0) \to P(\mu_0) \big) \Longrightarrow THH_*(W(k), (p)).
\]
We have essentially the same differentials as before, now with
\[
 d_{j+1}(x \mu_0^{p^j}) = x^{j+2} \mu_0^{p^j-1} \sigma x,
\]
and this tells us the following.

\begin{thm}
We have
\[
 THH_q(W(k), (p)) \cong \begin{cases} pW(k) & \textnormal{if $q=0$} \\ W_{\nu_p(i)+1}(k) & \textnormal{if $q=2i-1$ is odd} \\ 0 & \textnormal{if $q \geq 2$ is even} \end{cases}
\]
\end{thm}

In particular the long exact sequence coming from the fiber sequence defining $THH(W(k), (p))$ degenerates into short exact sequences
\begin{multline*} 0 \to THH_{2i}(k) \cong k \to THH_{2i-1}(W(k), (p)) \cong W_{\nu_p(i)+1}(k) \\ \to THH_{2i-1}(W(k)) \cong W_{\nu_p(i)}(k) \to 0. \end{multline*}
\end{example}

\begin{example} \label{ex:THHofW_nkrel}
Next we consider $THH(W_n(k), (p))$. Let
\[
 E_0^{*,*} = ker \big( P_n(x) \otimes E(\sigma x) \otimes \Gamma(x_n) \otimes P(\mu_0) \to P(\mu_0) \big)
\]
and let $d_0$ be generated multiplicatively by $d_0(\gamma_j(x_n)) = n x^{n-1} \gamma_{j-1}(x_n)$ for $j \geq 1$. Then we have a spectral sequence
\[
 E_1^{*,*} = H_*(E_0^{*,*}, d_0) \Longrightarrow THH_*(W_n(k), (p)).
\]
As long as $\nu_p(i)<n$ the following happens. The class $\mu_0^i$ was supposed to support a differential, but it is missing, so the target of the differential survives. This gives an extra class in $THH_{2i-1}(W_n(k), (p))$. If $\nu_p(i) \geq n$ then $\mu_0^i$ survives to give a class in $THH_{2i}(W_n(k))$; running the relative spectral sequence we then get one class less in $THH_{2i}(W_n(k), (p))$. Hence we find the following (compare Theorem \ref{t:THHofW_nk}).

\begin{thm} We have
\begin{eqnarray*}
THH_{2i}(W_n(k), (p)) & \cong & W_{\max(\nu_p(i), n-1)}(k) \oplus \bigoplus_{0 \leq j \leq i-1} W_{\max(\nu_p(j), n)}(k) \\
THH_{2i-1}(W_n(k), (p)) & \cong & W_{\max(\nu_p(i)+1, n)}(k) \oplus \bigoplus_{1 \leq j \leq i-1} W_{\max(\nu_p(j), n)}(k)
\end{eqnarray*}
\end{thm}
\end{example}

As in the DGA appearing in Theorem \ref{t:THHofW_nk} we will sometimes denote the class represented by $\sigma x$ in $THH_*(W(k), (p))$ or $THH_*(W_n(k), (p))$ by $\lambda_0$.

\section{The trace method} \label{s:tracemethod}
In this section we review the ``trace method'' for computing algebraic $K$-theory. Most of the material in this section is known, we include it here for the reader's convenience and for ease of reference. In some instances we have generalized known calculations from $\bF_p$ or $\bZ_p$ to $k$ or $W(k)$. We also introduce some notation that will be used in the later sections.

\subsection{Fixed points and geometric fixed points} \label{ss:fp}
From now on we will take $G$-spectrum to mean orthogonal $G$-spectrum in the sense of \cite{MaMa02}, see also \cite[Appendix B]{HHR}. Recall that an orthogonal $G$-spectrum is simply an orthogonal spectrum with a $G$-action; the difference between $G$-spectra indexed on different universes is taken care of by varying the model structure. In particular, different families of representations give different notions of fibrant replacement.

We will use a version of the positive complete stable model structure from \cite[Proposition B.63]{HHR} where the weak equivalences are defined using only finite subgroups of $G$. For details see \cite{ABGHLM}.

If we take the $H$-fixed points for some $H \leq G$ we get a $W(H)$-spectrum in the obvious way. It is important to note that taking fixed points does not commute with spectrification. In particular, if $X$ is a pointed $G$-space then $(\Sigma^\infty_G X)^H$ is very different from $\Sigma^\infty_{W(H)} X^H$. Instead, the classical tom Dieck splitting gives a formula for $(\Sigma^\infty_G X)^H$.

Second, we have the geometric fixed point spectrum $\Phi^H(E)$, obtained as a left Kan extension as in \cite[Definition V.4.3]{MaMa02}. If we apply geometric fixed points for a subgroup $H \leq G$ we again get a $W(H)$-spectrum. Taking geometric fixed points has the property that if $X$ is a pointed $G$-space then $(\Sigma^\infty_G X)^H \cong \Sigma^\infty_{W(H)} X^H$.

Now let $G=S^1$ and let $H = C_n$. If $E$ is an $S^1$-spectrum then $\Phi^{C_n}(E)$ is an $S^1/C_n$-spectrum. The $n$'th root provides an isomorphism $\rho_n : S^1 \to S^1/C_n$, and we can use this to change $\Phi^{C_n}(E)$ back into an $S^1$-spectrum $\rho_n^* \Phi^{C_n}(E)$.

The functor $\Phi^H$ preserves cofibrations and weak equivalences between cofibrant objects in our model structure, and has a left adjoint we will denote by $L\Phi^H$.

\begin{defn} [{\cite{BlMa}}]
A cyclotomic structure on an $S^1$-spectrum $E$ is a sequence of compatible maps
\[
 \rho_n^* \Phi^{C_n}(E) \to E
\]
for all $n \geq 2$ such that the induced map $\rho_n^* L\Phi^{C_n}(E) \to E$ is a weak equivalence. Here compatible means that the two maps from $\rho_{mn}^* \Phi^{C_{mn}}(E) = \rho_m^* \Phi^{C_m}(\rho_n^* \Phi^{C_n}(E))$ to $E$ agree.
\end{defn}

The canonical example of a cyclotomic spectrum is $\Sigma^\infty_{S^1} LX_+$, the equivariant suspension spectrum of a free loop space. In this case
\[
 \Phi^{C_n}(\Sigma^\infty_{S^1} LX_+) \simeq \Sigma^\infty_{S^1/C_n} (LX)^{C_n}_+,
\]
and we see that this is a cyclotomic spectrum because $(LX)^{C_n} \cong LX$.

We also know \cite{Bo1, HeMa97} that $THH(A)$ as defined by B\"okstedt is a cyclotomic spectrum; the proof that $THH(A)$ as defined in Section \ref{s:THH} above is a cyclotomic spectrum is similar to the proof found in \cite[Section 2]{HeMa97} in the classical setting. This should not be surprising, since
\[
 THH(\Sigma^\infty \Omega X_+) \simeq \Sigma^\infty LX_+.
\]

\begin{defn}
Let $A$ be a ring or symmetric ring spectrum. Then the $\TR$-groups of $A$ are the homotopy groups of the spectra
\[
 \TR^m(A) = THH(A)^{C_{p^{m-1}}}.
\]
(This is often denoted $\TR^m(A;p)$; we leave out the $p$ to simplify the notation.)
\end{defn}

These spectra are related by a number of maps, in a way that we now recall. There is a map $F : \TR^{m+1}(A) \to \TR^m(A)$ called Frobenius, which is given by inclusion of fixed points.

\begin{defn}
Let $A$ be a ring or symmetric ring spectrum. Then $\TF(A)$ is defined as
\[
 \TF(A) = \holim_F \TR^m(A).
\]

\end{defn}
The Frobenius has an associated transfer map $V : \TR^m(A) \to \TR^{m+1}(A)$ called the verschiebung. There is a map
\[
 d : \TR^m_q(A) \to \TR^m_{q+1}(A)
\]
defined as multiplication by the fundamental class of $S^1/C_{p^{m-1}}$.

Finally, there is a restriction map
\[
 R : \TR^{m+1}(A) \to \TR^m(A),
\]
which is defined using the cyclotomic structure on $THH(A)$. To be precise, the map
\[
 R : \TR^2(A) \to \TR^1(A)=THH(A)
\]
of non-equivariant spectra is given by the canonical map from fixed points to geometric fixed points, followed by the equivalence of the geometric fixed points with $THH(A)$. More generally $R : \TR^{m+1}(A) \to \TR^m(A)$ is the $C_{p^{m-1}}$ fixed points of this map. If we beef this up to include (virtual) $S^1$-representations the map $R$ takes the form
\[
 R : \Sigma^{\alpha} \TR^{m+1}(A) \to \Sigma^{\alpha'} \TR^m(A),
\]
where $\alpha=[\beta]-[\gamma] \in RO(S^1)$ and $\alpha'=\rho_p^* (\alpha^{C_p})$, see \cite{HeMa97b, Ger07}.

It is generally hard to understand fixed point spectra directly, and it is sometimes useful to compare the actual fixed point spectrum $\TR^{m+1}(A)$ to the homotopy fixed point spectrum $THH(A)^{hC_{p^m}}$. Let $T=THH(A)$, let $T_{hC_{p^m}}$ denote the homotopy orbit spectrum and let $T^{tC_{p^m}}$ denote the Tate spectrum. Then there is a fundamental diagram \cite[Theorem 1.10 and Section 2]{BoMa94}, as follows.

\[ \xymatrix{
 T_{hC_{p^m}} \ar[r]^-N \ar[d]^= & \TR^{m+1}(A) \ar[r]^-R \ar[d]^{\Gamma_m} & \TR^m(A) \ar[r]^-\partial \ar[d]^{\widehat{\Gamma}_m} & \Sigma T_{hC_{p^m}} \ar[d]^= \\
 T_{hC_{p^m}} \ar[r]^-{N^h} & T^{hC_{p^m}} \ar[r]^-{R^h} & T^{tC_{p^m}} \ar[r]^-\partial & \Sigma T_{hC_{p^m}}
} \]

If we take the homotopy inverse limit over $F$ we obtain a version of the fundamental diagram featuring $S^1$.

\[ \xymatrix{
 \Sigma T_{hS^1} \ar[r]^-N \ar[d]^= & \TF(A) \ar[d]^{\Gamma} \ar[r]^-R & \TF(A) \ar[r]^-\partial \ar[d]^{\widehat{\Gamma}} & \Sigma^2 T_{hS^1} \ar[d]^= \\
 \Sigma T_{hS^1} \ar[r]^-{N^h} & T^{hS^1} \ar[r]^-{R^h} & T^{tS^1} \ar[r]^-\partial & \Sigma^2 T_{hS^1}
} \]

Now consider the special case $A=\Sigma^\infty \Omega X_+$, so $THH(A) \simeq \Sigma^\infty_{S^1} LX_+$. The tom Dieck splitting says that
\[
 (\Sigma^\infty_{S^1} LX_+)^{C_{p^m}} \simeq \bigvee_{0 \leq k \leq m} (\Sigma^\infty LX_+)_{hC_{p^k}}.
\]
In this case the top row in the fundamental diagram splits. In general, the existence of the top row in the fundamental diagram can be thought of as a non-split version of the tom Dieck splitting for general $A$.

Finally we get to topological cyclic homology.

\begin{defn}
Let $A$ be a ring or symmetric ring spectrum. The topological cyclic homology $\TC(A)$ of $A$ is the homotopy equalizer
\[
 \TC(A) \to \TF(A) \underset{id}{\overset{R}{\rightrightarrows}} \TF(A).
\]
Alternatively, it can be defined as the homotopy equalizer
\[
 \TC(A) \to \TR(A) \underset{id}{\overset{F}{\rightrightarrows}} \TR(A),
\]
where $\TR(A) = \holim_R \TR^m(A)$, or as $\TC(A) = \holim_{R, F} \TR^m(A)$.
\end{defn}

There is a trace map
\[
 trc : K(A) \to \TC(A)
\]
which is an isomorphism on homotopy groups in degree $\geq 0$ after $p$-completion if $A$ is e.g.\ a finite $W(k)$-algebra \cite{Mc97}. These comparison results go through relative $\TC$ and relative $K$-theory.

Given a functor $F$ from rings (or symmetric ring spectra) to spectra and an ideal $I \subset A$, we define $F(A,I)$ as the homotopy fiber
\[
 F(A,I) \to F(A) \to F(A/I).
\]
This defines relative $K$-theory and $\TC$, and we have a relative trace map
\[
 trc : K(A,I) \to \TC(A,I).
\]
What McCarthy \cite{Mc97} actually shows is that this relative trace map is an equivalence after $p$-completion when $I$ is nilpotent. (Actually the relative trace map is an equivalence even before $p$-completing, see \cite{DGM} for details.)

The calculation of $\TC(k)$ recalled below plus Kratzer's calculation of $K(k)$ \cite{Kr80} provides the base case which we use to conclude that the absolute trace map is an equivalence in non-negative degrees after $p$-completion for certain rings.

In particular this means that up to $p$-completion we have
\[
 K_q(W_n(k), (p)) \cong \TC_q(W_n(k), (p))
\]
for all $q$.

In some cases we can use a result of Tsalidis to study $\TR^m(A)$ in terms of the $C_{p^m}$ Tate spectrum.

\begin{thm}[Tsalidis, \cite{Ts98}] \label{t:Tsalidis}
Let $A$ be a connective symmetric ring spectrum of finite type. Suppose
\[
 \wh{\Gamma}_1 : \pi_q THH(A) \to \pi_q THH(A)^{tC_p}
\]
is an isomorphism for $q \geq q_0$. Then
\[
 \wh{\Gamma}_m : \TR^m_q(A) \to \pi_q THH(A)^{tC_{p^m}}
\]
is an isomorphism for $q \geq q_0$ for all $m$.
\end{thm}

This allows for an induction argument, as follows. Recall \cite{GrMa95, BoMa94} that there is a Tate spectral sequence converging to $\pi_* THH(A)^{tC_{p^m}}$, and that we get spectral sequences converging to $\pi_* THH(A)_{hC_{p^m}}$ and $\pi_* THH(A)^{hC_{p^m}}$ by (with a small modification in filtration $0$) restricting to the first or second quadrant, respectively. If the conditions of Tsalidis' Theorem hold and we understand $\TR^m_*(A)$, we can often understand the spectral sequence converging to $\pi_* THH(A)^{tC_{p^m}}$ because we know what it converges to in degree $q \geq q_0$. Then restricting this spectral sequence to the second quadrant gives a spectral sequence computing $\pi_* THH(A)^{hC_{p^m}}$, and this determines $\TR^{m+1}_q(A)$ for $q \geq q_0$.

By taking the homotopy inverse limit over $F$, we can also conclude that the maps $\Gamma : \TF_q(A) \to \pi_q THH(A)^{hS^1}$ and $\widehat{\Gamma} : \TF_q(A) \to \pi_q THH(A)^{tS^1}$ are isomorphisms for $q \geq q_0+1$.

\subsection{Topological cyclic homology of $k$} \label{ss:TCk}
Many computations rely on the corresponding computations for $k$, so following \cite{HeMa97} we spell this case out first. Recall that $THH_*(k)=P(\mu_0)$ is a polynomial algebra over the ground field $k$ on a degree $2$ generator $\mu_0$. Then the Tate spectral sequence converging to $\pi_* THH(k)^{tC_{p^m}}$ takes the following form:
\[
 \widehat{E}_2^{*,*} = P(\mu_0) \otimes E(u_m) \otimes P(t,t^{-1}) \Longrightarrow \pi_* THH(k)^{tC_{p^m}}.
\]
This is bigraded by fiber degree and homological degree, with $|\mu_0|=(2,0)$, $|u_m|=(0,-1)$ and $|t|=(0,-2)$. The topological degree is the sum of the two degrees. The class $v_0=t\mu_0$ represents multiplication by $p$ and is a permanent cycle. We have a differential
\[
 d_{2m+1}(u_m) \doteq t^{m+1} \mu_0^m = t v_0^m,
\]
leaving
\[
 \widehat{E}_{2m+2}^{*,*}=\widehat{E}_\infty^{*,*} = P_m(v_0) \otimes P(t,t^{-1}).
\]
This is the associated graded of
\[
 \pi_* THH(k)^{tC_{p^m}} \cong W_m(k)[t,t^{-1}].
\]
When $m=1$ the map $\widehat{\Gamma}_1 : THH_*(k) \to \pi_* THH(k)^{tC_p}$ is given by $\wh{\Gamma}(\mu_0) \doteq t^{-1}$. This is an isomorphism in non-negative degrees and Tsalidis' Theorem applies.

To compute $\pi_* THH(k)^{hC_{p^m}}$ we restrict the Tate spectral sequence to the second quadrant, and we have
\[
 E_2^{*,*} = P(\mu_0) \otimes E(u_m) \otimes P(t).
\]
We have the same $d_{2m+1}$-differential, which leaves
\[
 E_{2m+2}^{*,*}=E_\infty^{*,*} = P_m(v_0)\{t^i \,\, | \,\, i>0\} \oplus P_{m+1}(v_0)\{\mu_0^j \,\, | \,\, j \geq 0\}.
\]
This is the associated graded of
\[
 \pi_* THH(k)^{hC_{p^m}} \cong W_m(k)\{t^i \,\, | \,\, i>0\} \oplus W_{m+1}(k)\{\mu_0^j \,\,  | \,\, j \geq 0\}.
\]

Next we compute $R : \pi_* THH(k)^{C_{p^m}} \to \pi_* THH(k)^{C_{p^{m-1}}}$, and here we need to be a little bit careful. From \cite{HeMa97} we know that we have an isomorphism
\[
 \rho^R_m : \pi_0 THH(k)^{C_{p^{m-1}}} \to W_m(k)
\]
which is compatible with the restriction map $R$. But we would like to take the inverse limit over $F$ rather than $R$, because we want to study $\TF(A)$ rather than $\TR(A)$.

We have a commutative diagram
\[ \xymatrix{
 \ldots \ar[r] & W_3(k) \ar[r]^-F \ar[d]^{W_3(\phi^3)} & W_2(k) \ar[r]^-F \ar[d]^{W_2(\phi^2)} & k \ar[d]^\phi \\
 \ldots \ar[r] & W_3(k) \ar[r]^-R & W_2(k) \ar[r]^-R & k
} \]
where $W_m(\phi^m) : W_m(k) \to W_m(k)$ is an isomorphism, and for our purposes it is better to use the isomorphism
\[
 \rho^F_m : \pi_0 THH(k)^{C_{p^{m-1}}} \to W_m(k)
\]
given by $\rho^F_m = W_m(\phi^m) \circ \rho^R_m$. This extends to an isomorphism
\[
 \rho^F_m : \pi_* THH(k)^{C_{p^{m-1}}} \to W_m(k)[\mu_0].
\]

\begin{lemma} \label{l:RandFconventions}
Suppose we use the above map $\rho^F_m$ to identify $\pi_* THH(k)^{C_{p^{m-1}}}$ with $W_m(k)[\mu_0]$. Then the Frobenius map $F : THH(k)^{C_{p^m}} \to THH(k)^{C_{p^{m-1}}}$ is given by the usual restriction map on Witt vectors and by $\mu_0 \mapsto \mu_0$. The restriction map $R : THH(k)^{C_{p^m}} \to THH(k)^{C_{p^{m-1}}}$ is given by the usual restriction map on Witt vector followed by $W_m(\phi^{-1})$ and by $\mu_0 \mapsto p \lambda_m \mu_0$ for some unit $\lambda_m \in \bZ/p^m$.
\end{lemma}

It follows that $\TF_*(k) \cong W(k)[\mu_0]$ and that
\[
 \TC_i(k) \cong \begin{cases} \bZ_p & \textnormal{if $i=0$} \\ coker(W(\phi)-1) & \textnormal{if $i = -1$} \\ 0 & \textnormal{otherwise} \end{cases}
\]
Here we have used that the (co)equalizer of $W(\phi^{-1})$ and $1$ is isomorphic to the (co)kernel of $W(\phi) - 1$.

The trace map $K_*(k) \to \TC_*(k)$ is, after $p$-completion, an isomorphism in degree $0$ and trivial in degree $-1$, since $K(A)$ is a connective spectrum for any ring $A$.

Together with Kratzer's calculation \cite[Corollary 5.5]{Kr80} of $K(k)$ this provides the base case where the trace map is an equivalence on non-negative homotopy groups after $p$-adic completion.

\subsection{Topological cyclic homology of $W(k)$}
Next we turn to the topological cyclic homology of $W(k)$. In this case the Tate spectral sequence converging to $\pi_* THH(W(k))^{tS^1}$ has $E_2$-term
\[
 \wh{E}_2^{*,*} = \Big( W(k)\{1\} \oplus \bigoplus_{i \geq 1} W_{\nu_p(i)+1}(k)\{\lambda_1\mu_1^{i-1}\} \Big) \otimes P(t,t^{-1}).
\]
From \cite{BoMa94} we know the behavior of the corresponding spectral sequence with mod $p$ coefficients, and this lets us say some things about this spectral sequence. But rather than going into details we will instead study the relative version, which has the form
\[
 \wh{E}_2^{*,*} = \Big( W(k)\{p\} \oplus \bigoplus_{i \geq 1} W_{\nu_p(i)+1}(k)\{\lambda_0 \mu_0^{i-1}\} \Big) \otimes P(t,t^{-1}).
\]

\begin{prop} \label{p:TateSSforWk}
The differentials in the above spectral sequence are all $P(t,t^{-1})$-linear and given by
\[
 d_{2r}(p^r) \doteq r t^r \lambda_0 \mu_0^{r-1}.
\]
The $E_\infty$ term is given by
\[
 \wh{E}_\infty^{*,*} = \bigoplus_{i \geq 1} W_{\nu_p(i)+1}(k)\{\lambda_1 \mu_1^{i-1}\} \otimes P(t,t^{-1}).
\]

\end{prop}

\begin{proof}
We have a long exact sequence of spectral sequences
\[
 \ldots \to \wh{E}_r^{*,*}(W(k), (p)) \to \wh{E}_r^{*,*}(W(k)) \to \wh{E}_r^{*,*}(k) \to \ldots
\]
for each $r \geq 2$. Consider the class $p^r \in \wh{E}_2^{0,0}(W(k))$. In the Tate spectral sequence for $k$ there is a hidden extension, with $p^r$ represented by $t^r \mu_0^r$. The connecting map $THH_*(k) \to THH_{*-1}(W(k), (p))$ maps $\mu_0^r$ to $r \lambda_0 \mu_0^{r-1}$, hence in the Tate spectral sequence it maps $t^r \mu_0^r$ to $r t^r \lambda_0 \mu_0^{r-1}$. If there was no such $d_{2r}$ differential this would break exactness at the $E_r$ term, so the result follows.
\end{proof}

\begin{remark}
The above result can be explained in a simple but non-rigorous way as follows. We have $\pi_* THH(k)^{tS^1} \cong W(k)[t,t^{-1}]$, and by ``compressing'' the filtration in $\wh{E}_2^{*,*}(k)$ the map from $\wh{E}_2^{*,*}(W(k))$ becomes surjective with kernel concentrated in odd total degree, isomorphic to the given $E_\infty$ term.
\end{remark}

This computes $\TF_*(W(k), (p))$ up to extensions in two ways, using either the map $\Gamma : \TF_*(W(k), (p)) \to \pi_* THH(W(k), (p))^{hS^1}$ or the map $\wh{\Gamma} : \TF_*(W(k), (p)) \to \pi_* THH(W(k), (p))^{tS^1}$. One could now attempt to compute $\TC_*(W(k), (p))$ by understanding the restriction map $R$ on $\TF_*(W(k), (p))$. This is complicated, and we do not now how to do this without passing to mod $p$ coefficients and following \cite{BoMa94, BoMa95}. We will omit a discussion of the mod $p$ calculation as we will not need it.

\subsection{Topological cyclic homology of $k[x]/(x^n)$} \label{ss:TCtruncated}
We also need the computation of $\TC_*(k[x]/(x^n))$ from \cite{HeMa97b}. Suppose $\Pi$ is a pointed monoid, and let $k(\Pi)$ denote the pointed monoid algebra. Then $THH(k(\Pi)) \simeq THH(k) \sma B^{cy}_\sma(\Pi)$, and this is an equivalence of $S^1$-equivariant spectra. In particular, let $\Pi_n=\{0,1,x,\ldots,x^{n-1}\}$ so that $k(\Pi_n)=k[x]/(x^n)$. Then it is clear that $B^{cy}_\sma(\Pi_n)$ splits as a wedge of homogeneous summands, using the degree in $x$, and Hesselholt and Madsen calculated the $S^1$-equivariant homotopy type of as follows:

\begin{thm}[Hesselholt-Madsen \cite{HeMa97b}] \label{t:HMhomotopytype}
The cyclic bar construction $B^{cy}(\Pi_n)$ splits, $S^1$-equivariantly, as
\[
 B^{cy}(\Pi_n) \cong \bigvee_{s \geq 0} B^{cy}(\Pi_n;s),
\]
where $B^{cy}(\Pi_n;0) = S^0$,
\[
 B^{cy}(\Pi_n;s) \simeq S^1(s)_+ \sma S^{\lambda_d}
\]
if $n$ does not divide $s$ and $B^{cy}(\Pi_n;s)$ sits as the mapping cone
\[
 S^1(s/n)_+ \sma S^{\lambda_d} \overset{n}{\to} S^1(s)_+ \sma S^{\lambda_d} \to B^{cy}(\Pi_n;s)
\]
if $n$ divides $s$.
\end{thm}

Here $d=\lfloor \frac{s-1}{n} \rfloor$, $\lambda_d = \bC(1) \oplus \ldots \oplus \bC(d)$, and $S^1(s)$ denotes $S^1$ as an $S^1$-space with an accelerated action. If $p$ does not divide $n$ then this simplifies after $p$-completion as
\[
 B^{cy}(\Pi_n)^\wedge_p \simeq (S^0)^\wedge_p \vee \bigvee_{n \nmid s} B^{cy}(\Pi_n;s)^\wedge_p.
\]
The homotopy groups of $THH(k[x]/(x^n))$ are given by the homology of the DGA
\[
 P_n(x) \otimes E(\sigma x) \otimes \Gamma(x_n) \otimes P(\mu_0),
\]
with differential generated multiplicatively by $d(x_n)=nx^{n-1} \sigma x$. (Compare with Theorem \ref{t:THHofW_nk}.)

Because the above splitting is $S^1$-equivariant it follows that
\[
 \TR^m(k[x]/(x^n)) \simeq \bigvee_{s \geq 0} \TR^m(k[x]/(x^n); s)
\]
and that
\[
 \TF(k[x]/(x^n)) \simeq \bigvee_{s \geq 0} \TF(k[x]/(x^n);s).
\]

We need to understand the Tate spectrum and homotopy fixed point spectrum of $THH(k[x]/(x^n);s)$. We start with the Tate spectrum. It is convenient to define
\[
 z_s = t^{\lfloor s/n \rfloor} \gamma_{\lfloor s/n \rfloor}(x_n) x^{\{s/n\}},
\]
where $\{s/n\}$ denotes the residue of $s$ mod $n$. This class lives in total degree $0$ in the Tate spectral sequence. Then we can describe the $E_2$ term of the Tate spectral sequence converging to $\pi_* THH(k[x]/(x^n);s)^{tS^1}$ as follows. We start with
\[
 \wh{E}_0^{*,*}(s) = P(\mu_0) \otimes P(t,t^{-1}) \otimes k\{z_s, z_{s-1} \sigma x\}.
\]
If $n \mid s$ and $p \nmid n$ we define $d_0(z_s) = tz_{s-1} \sigma x$, otherwise we define $d_0=0$. Then $\wh{E}_2^{*,*}(s)$ is the homology of $(\wh{E}_0^{*,*}(s), d_0)$. If $n \mid s$ and $p \nmid n$ we have $\wh{E}_2^{*,*}(s)=0$.

The higher differentials are given by
\[
 d_{2\nu_p(s)+2}(z_s) \doteq t^{\nu_p(s)+1} \mu_0^{\nu_p(s)} z_{s-1} \sigma x
\]
if $n \nmid s$ and
\[
 d_{2\nu_p(n)}(z_s) \doteq t^{\nu_p(n)+1} \mu_0^{\nu_p(n)} z_{s-1} \sigma x
\]
if $n \mid s$. (This includes the above $d_0$ differential in the case $\nu_p(n)=0$.)

This leaves the $E_\infty$ term
\[
 \wh{E}_\infty^{*,*}(s) = P_{\nu_p(s)}(v_0) \otimes P(t,t^{-1}) \otimes k\{z_{s-1} \sigma x\} 
\]
if $n \nmid s$ and
\[
 \wh{E}_\infty^{*,*}(s) = P_{\nu_p(n)}(v_0) \otimes P(t,t^{-1}) \otimes k\{z_{s-1} \sigma x\}
\]
if $n \mid s$. This is the associated graded of
\[
 \pi_* THH(k[x]/(x^n);s)^{tS^1} = W_{\nu_p(s)}(k) \otimes P(t,t^{-1}) \otimes k\{z_{s-1} \sigma x\}
\]
if $n \nmid s$ and
\[
 \pi_* THH(k[x]/(x^n);s)^{tS^1} = W_{\nu_p(n)}(k) \otimes P(t,t^{-1}) \otimes k\{z_{s-1} \sigma x\}
\]
if $n \mid s$.

Restricting the Tate spectral sequence to the second quadrant we obtain a calculation of $\pi_* THH(k[x]/(x^n);s)^{hS^1}$. It will be convenient to introduce another family of elements. Let
\[
 y_s = \mu_0^{-\lfloor s/n \rfloor} \gamma_{\lfloor s/n \rfloor}(x_n) x^{\{s/n\}}.
\]
Then the class $\mu_0^j y_s$ is in $THH_{2j}(k[x]/(x^n);s)$ for $j \geq \lfloor \frac{s}{n} \rfloor$ and the class $\mu_0^j y_{s-1} \sigma x$ is in $THH_{2j+1}(k[x]/(x^n);s)$ for $j \geq d$. (Recall that we defined $d = \lfloor \frac{s-1}{n} \rfloor$.) With this naming convention we then get
\begin{eqnarray*}
 \pi_* THH(k[x]/(x^n);s)^{hS^1} = & \bigoplus_{i > 0} & W_{\nu_p(s)}(k)\{t^i \mu_0^d y_{s-1} \sigma x\} \\
 \oplus & \bigoplus_{j \geq d} & W_{\nu_p(s)+1}(k)\{\mu_0^j y_{s-1} \sigma x\}
\end{eqnarray*}
if $n \nmid s$ and
\begin{eqnarray*}
 \pi_* THH(k[x]/(x^n);s)^{hS^1} = & \bigoplus_{i > 0} & W_{\nu_p(n)}(k)\{t^i \mu_0^d y_{s-1} \sigma x\} \\
 \oplus & \bigoplus_{j \geq d} & W_{\nu_p(n)}(k)\{\mu_0^j y_{s-1} \sigma x\}
\end{eqnarray*}
if $n \mid s$.

The above calculation of $\pi_* THH(k[x]/(x^n))^{hS^1}$ and $\pi_* THH(k[x]/(x^n))^{tS^1}$ does not, in itself, compute $\TF_*(k[x]/(x^n))$, because Tsalidis' Theorem does not apply. But it is possible to compute $\TF_*(k[x]/(x^n); s)$ for each $s$ directly, identifying it with $\TR^{\nu_p(s)+1}_{*-\lambda_d-1}(k)$ if $n \nmid s$ and with the cokernel of $V^{\nu_p(n)} : \TR^{\nu_p(s/n)+1}_{*-\lambda_d-1}(k) \to \TR^{\nu_p(s)+1}_{*-\lambda_d-1}(k)$ if $n \mid s$. And we have the following computation, see \cite{HeMa97b}. See also \cite{Ger07, AnGe11} in the case $k=\bF_p$.

\begin{thm}
Let $\lambda$ be an actual complex $S^1$-representation. Then $\TR^m_{*-\lambda}(k)$ is concentrated in even degree. If $i \geq \dim_\bC(\lambda)$ we have $\TR^m_{2i-\lambda}(k)=W_m(k)$. If $\dim_\bC(\lambda^{(j-1)}) > i \geq \dim_\bC(\lambda^{(j)})$ then $\TR^n_{2i-\lambda}(k)=W_{m-j}(k)$.
\end{thm}

This is proved using an $RO(S^1)$-graded version of the fundamental diagram. For any virtual $S^1$-representation $\alpha$ we have a fundamental diagram
\[ \xymatrix{
 T[\alpha]_{hC_{p^m}} \ar[r]^-N \ar[d]^= & \TR^{m+1}(A)[\alpha] \ar[r]^-R \ar[d]^{\Gamma_m} & \TR^m(A)[\alpha'] \ar[d]^{\widehat{\Gamma}_m} \\ T[\alpha]_{hC_{p^m}} \ar[r]^-{N^h} & T[\alpha]^{hC_{p^m}} \ar[r]^-{R^h} & T[\alpha]^{tC_{p^m}}
} \]
This diagram can also be used to compute $R : \TR^{m+1}_{*-\lambda}(k) \to \TR^m_{*-\lambda'}(k)$.

We use Theorem \ref{t:HMhomotopytype} above and find (compare \cite[Section 8.2]{HeMa97}) that if $n \nmid s$ then
\begin{multline*}
\TF(A[x]/(x^n);s) \simeq (S^1(s)_+ \sma S^{\lambda_d} \sma THH(A))^{S^1} \\ \simeq \Sigma F(S^1(s)_+, THH(A) \sma S^{\lambda_d})^{S^1} \simeq \Sigma (THH(A) \sma S^{\lambda_d})^{C_s}
\end{multline*}
up to $p$-completion. Similarly, if $n \mid s$ then $\TF(A[x]/(x^n);s)$ sits in a cofibration sequence
\[
 \Sigma (THH(A) \sma S^{\lambda_d})^{C_{s/n}} \overset{V_n}{\longrightarrow} \Sigma (THH(A) \sma S^{\lambda_d})^{C_s} \to \TF(A[x]/(x^n);s).
\]
Hence
\[
 \TF_*(k[x]/(x^n); s) \cong \TR^{\nu_p(s)+1}_{*-1-\lambda_d}(k)
\]
when $n \nmid s$ and similarly for the case $n \mid s$. This is what Hesselholt and Madsen used to compute $K_*(k[x]/(x^n))$.

With this we can describe the maps $\Gamma : \TF_*(k[x]/(x^n)) \to THH(k[x]/(x^n))^{hS^1}$ and $\widehat{\Gamma} : \TF_*(k[x]/(x^n)) \to THH(k[x]/(x^n))^{tS^1}$. The map $\Gamma$ sends $\TF(k[x]x/(x^n); s)$ to $THH(k[x]/(x^n); s)^{hS^1}$ and even though Tsalidis' theorem does not apply we do have the following:

\begin{thm} \label{t:Gamma_injective}
In degree $2i+1$ for $i \geq d$ the map
\[
 \Gamma : \TF_{2i+1}(k[x]/(x^n); s) \to \pi_{2i+1} THH(k[x]/(x^n); s)^{hS^1}
\]
is an isomorphism. In degree $2i+1$ for $i < d$ the map
\[
 \Gamma : \TF_{2i+1}(k[x]/(x^n); s) \to \pi_{2i+1} THH(k[x]/(x^n); s)^{hS^1}
\]
is injective.
\end{thm}

We have a similar result for the map $\widehat{\Gamma}$. In this case $\widehat{\Gamma}$ sends $\TF(k[x]/(x^n); s)$ to $THH(k[x]/(x^n); ps)$.

\begin{thm} \label{t:Gammahat_injective}
In degree $2i+1$ for $i \geq d$ the map
\[
 \widehat{\Gamma} : \TF_{2i+1}(k[x]/(x^n); s) \to \pi_{2i+1} THH(k[x]/(x^n); ps)^{tS^1}
\]
is an isomorphism. In degree $2i+1$ for $i<d$ the map
\[
 \widehat{\Gamma} : \TF_{2i+1}(k[x]/(x^n); s) \to \pi_{2i+1} THH(k[x]/(x^n); ps)^{tS^1}
\]
is injective. 
\end{thm}

If we use Theorem \ref{t:Gamma_injective} to name elements in $\TF_*(k[x]/(x^n))$, the map $\wh{\Gamma}$ is given by
\[
 \wh{\Gamma}(\mu_0^a y_{s-1} \sigma x) = t^{-a} z_{ps-1} \sigma x.
\]

From this we can read off the action of
\[
 R : \TF_{2i+1}(k[x]/(x^n); s) \to \TF_{2i+1}(k[x]/(x^n); s/p).
\]

\begin{thm} \label{t:descriptionofR}
Suppose $\nu_p(s) \geq 1$. In degree $2i+1$ for $i \geq d$ the map
\[
 R : \TF_{2i+1}(k[x]/(x^n); s) \to \TF_{2i+1}(k[x]/(x^n); s/p)
\]
is multiplication by $p^{i-d}$. In degree $2i+1$ for $i<d$ the map $R$ is an isomorphism.
\end{thm}

In particular this means that there is a stable range. If $i < d$ then
\[
 R : \TF_{2i+1}(k[x]/(x^n); s) \to \TF_{2i+1}(k[x]/(x^n); s/p)
\]
is an isomorphism, and if $i = d$ it is surjective.

\section{A spectral sequence on fixed points}
In this section we introduce the filtered ring spectral sequences on fixed points and use it to calculate the associated graded of $\TF_*(W_n(k), (p))$. Using this we prove Theorem \ref{t:mainGalois}.

\subsection{A spectral sequence for $\TR^n(A)$}
The spectral sequence in Theorem \ref{t:THHSS} comes from an $S^1$-equivariant filtration on $THH(A)$, so it is reasonable to expect it to induce a filtration on fixed points as well. Once we have this, we get an induced spectral sequence on fixed points as well.

\begin{thm} \label{t:TRSS}
Suppose $A$ is a complete filtered ring or symmetric ring spectrum with associated graded $Gr A$. Then there is a weakly convergent spectral sequence
\[
 E_1^{s,t} = \TR^m_{s+t}(Gr A; s) \Longrightarrow \TR^m_{s+t}(A).
\]
If $A$ is commutative then this is an algebra spectral sequence.
\end{thm}

\begin{proof}
We prove the case $m=2$, the general case is similar. We use the $p$-fold edgewise subdivision model of $THH$, which is the $S^1$-spectrum with $V$'th space the geometric realization of
\[
 THH^{[p]}(A)_q(V) = \underset{I^{p(q+1)}}{\hocolim} \, \Omega^{i_0 +\ldots+i_{p(q+1)-1}} \big( A_{i_0} \sma \ldots \sma A_{i_{p(q+1)-1}} \sma S^V \big).
\]
The advantage of this model is that we have a simplicial action of $C_p$.

While $THH(A)$ might not be fibrant, the discussion in \cite[Section 2 and Appendix A]{HeMa97} implies that the fixed points of $THH(A)$ (without fibrant replacement) calculates the correct homotopy type.

We have a filtration on each $THH^{[p]}(A)_q(V)$ coming from the filtration on each space $A_i$ in the spectrum $A$, and this induces a filtration on $THH^{[p]}(A)$ which is equivalent to the filtration on $THH(A)$ considered before. With this model it is clear that taking fixed points preserves the filtration, since the representation spheres $S^V$ are all in filtration $0$.
\end{proof}

There is of course a similar spectral sequence converging to the homotopy groups of the relative spectrum.

\begin{cor}
Suppose $A$ is a complete filtered ring or symmetric ring spectrum with associated graded $Gr A$ and let $I=F^1 A \subset A$. Then there is a spectral sequences
\[
 E_1^{s,t} = \begin{cases} \TR^m_{s+t}(Gr A; s) & \textnormal{if $s \geq 1$} \\ 0 & \textnormal{if $s=0$} \end{cases} \Longrightarrow \TR^m_{s+t}(A, I).
\]
\end{cor}

A description of the $E_1$-term of this spectral sequence for $(A,I)=(W_n(k), (p))$ follows from the calculations in \cite{HeMa97b}, recalled in Section \ref{ss:TCtruncated} above. Because we will only need the corresponding spectral sequence for $\TF$ we omit the details.

\subsection{A spectral sequence for $\TF(A)$}
The Frobenius $F$ is simply the inclusion of fixed points, so it is compatible with the filtration and we can take a homotopy inverse limit to get a spectral sequence converging to $\TF_*(A)$.

\begin{thm}
Suppose $A$ is a complete filtered ring or symmetric ring spectrum with associated graded $Gr A$. Then there is a weakly convergent spectral sequence
\[
 E_1^{s,t} = \TF_{s+t}(Gr A; s) \Longrightarrow \TF_{s+t}(A).
\]
\end{thm}

As usual there is a relative version.

\begin{cor}
Suppose $A$ is a complete filtered ring or symmetric ring spectrum with associated graded $Gr A$ and let $I=F^1 A \subset A$. Then there is a spectral sequences
\[
 E_1^{s,t} = \begin{cases} \TF_{s+t}(Gr A; s) & \textnormal{if $s \geq 1$} \\ 0 & \textnormal{if $s=0$} \end{cases} \Longrightarrow \TF_{s+t}(A, I).
\]
\end{cor}

For $A=W_n(k)$ this $E_1$-term is studied in \cite{HeMa97b} as recalled in Section \ref{ss:TCtruncated} above, and we find the following.

\begin{prop}
Suppose $p \nmid n$. Then the above spectral sequence converging to $\TF_*(W_n(k), (p))$ has $E_1$-term
\[
 E_1^{s,*} = \TR^{\nu_p(s)+1}_{*-\lambda_d-1}(k) \qquad \textnormal{if $n \nmid s$}
\]
and $E_1^{s,*}=0$ if $n \mid s$ for $s \geq 1$.
\end{prop}

Note that this is concentrated in odd topological degree, and hence this spectral sequence collapses at the $E_1$-term.
In particular, $E_1^{s,*}$ is a $W_{\nu_p(s)+1}(k)$ in sufficiently high odd total degree.

\begin{prop}
Suppose $p \mid n$. Then the above spectral sequence converging to $\TF_*(W_n(k), (p))$ has $E_1$-term
\[
 E_1^{s,*} = \begin{cases} \TR^{\nu_p(s)+1}_{*-\lambda_d-1}(k) & \qquad \textnormal{if $n \nmid s$} \\ coker \big( \TR^{\nu_p(s/n)+1}_{*-\lambda_d-1}(k) \overset{V^{\nu_p(n)}}{\longrightarrow} \TR^{\nu_p(s)+1}_{*-\lambda_d-1}(k) \big) & \qquad \textnormal{if $n \mid s$} \end{cases}
\]
for $s \geq 1$.
\end{prop}

In the case $n \mid s$ the cokernel is isomorphic to $W_{\nu_p(n)}(k)$ in sufficiently high odd total degree, and again we see that the $E_1$-term is concentrated in odd topological degree.

\begin{cor} \label{c:collapses}
The spectral sequence converging to $\TF_*(W_n(k), (p))$ collapses at the $E_1$-term.
\end{cor}

We compare this to $W(k)$, for which we find the following. (See also Proposition \ref{p:TateSSforWk} above.)

\begin{cor} \label{c:collapseWk}
The spectral sequence converging to $\TF_*(W(k), (p))$ has $E_1$-term
\[
 E_1^{s,*} = \TR^{\nu_p(s)+1}_{*-1}(k)
\]
for $s \geq 1$. This spectral sequence also collapses at the $E_1$-term.
\end{cor}

\begin{proof}[Proof of Theorem \ref{t:mainGalois}]
Suppose $k \to k'$ is a $G$-Galois extension of perfect fields of characteristic $p$ for a finite group $G$. Then it follows from Corollary \ref{c:collapses} and \ref{c:collapseWk} that
$\TF_*(W_n(k')) \cong \TF_*(W_n(k)) \otimes_{W(k)} W(k')$ with the induced $G$-action. Hence the homotopy fixed point spectral sequence
\[
 H^*(G; \TF_*(W_n(k'))) \Longrightarrow \pi_*(\TF(W_n(k'))^{hG})
\]
collapses at the $E_2$-term, and it follows that the canonical map
\[
 \TF(W_n(k)) \to \TF(W_n(k'))^{hG}
\]
is an equivalence.

The maps $R$ and $1$ are $G$-equivariant, and homotopy equalizers commute with homotopy fixed points. Hence the canonical map
\[
 \TC(W_n(k)) \to \TC(W_n(k'))^{hG}
\]
is an equivalence as well. The statement of the theorem follows by taking connective covers and $p$-completing.
\end{proof}

\section{A commuting square of spectral sequences} \label{s:commutingsquareofSSs}
In this section we compare two ways of calculating the homotopy groups of the Tate spectrum $THH(W_n(k), (p))^{tS^1}$. The main purpose is to better understand the map $\wh{\Gamma} : \TF_*(W_n(k), (p)) \to \pi_* THH(W_n(k), (p))^{tS^1}$ and thereby better understand the restriction map on $\TF_*(W_n(k), (p))$.

\subsection{A commuting square} \label{ss:commutingsquaregeneral}
We begin with a general observation about what happens when we have a ``commuting square'' of spectral sequences. We claim no originality for this, but because we could not find exactly what we need in the literature we include this material here. But see \cite{AnLi15} for a related discussion.

Suppose we have a spectrum $X$ with two compatible filtrations on it. That means we have a spectrum $F^{i,j} X$ for each $i$ and $j$, with maps $F^{i,j} X \to F^{i-1,j} X$ and $F^{i,j} X \to F^{i,j-1} X$ which commute in the obvious sense, with $X = F^{-\infty, -\infty} X$ the homotopy colimit. By forgetting one of the filtrations we get a horizontal filtration
\[
 \ldots \to F^{i+1,-\infty} X \to F^{i,-\infty} X \to F^{i-1,-\infty X} \to \ldots
\]
and a corresponding horizontal spectral sequence, and by forgetting the other we get a vertical filtration
\[
 \ldots \to F^{-\infty, j+1} X \to F^{-\infty, j} X \to F^{-\infty, j-1} \to \ldots
\]
and a corresponding vertical spectral sequence.

Similarly, for each $j$ we have a horizontal filtration
\[
 \ldots \to F^{i+1,j} X/F^{i+1,j+1} X \to F^{i,j} X/F^{i,j+1} X \to F^{i-1,j} X/F^{i-1,j+1} X \to \ldots
\]
and a corresponding horizontal spectral sequence, and for each $i$ we have a vertical filtration
\[
 \ldots \to F^{i,j+1} X/F^{i+1,j+1} X \to F^{i,j} X/F^{i+1,j} X \to F^{i,j-1} X/F^{i+1,j-1} X \to \ldots
\]
and a corresponding vertical spectral sequence.

Putting all of these together we get a ``commuting square'' of spectral sequences
\[ \xymatrix{
 E_1^{i,j,k} = \pi_{i+j+k} Gr^{i,j} X \ar@{=>}[r] \ar@{=>}[d] & (E_1')^{i+k,j} = \pi_{i+j+k} Gr^j_v X \ar@{=>}[d] \\
 (E_1'')^{i,j+k} = \pi_{i+j+k} Gr^i_h X \ar@{=>}[r] & \pi_{i+j+k} X.
} \]
Here $Gr^{i,j} X$ is the iterated homotopy cofiber $\frac{F^{i,j} X/F^{i+1,j} X}{F^{i+1,j} X/F^{i+1,j+1} X}$, while $Gr^i_h X$ is the homotopy cofiber $F^{i,-\infty}/F^{i+1,-\infty}$ and similarly for $Gr^j_v X$.

\begin{question}
Given a class in $\alpha \in \pi_m X$, how can we find all the representatives of $\alpha$ in $E_1^{*,*,*}$?
\end{question}

An element $y \in E_1^{i,j,m-i-j}$ represents $\alpha$ if $y$ lifts to an element $\tilde{y} \in \pi_m F^{i,j} X$ which maps to $\alpha$ under the map $F^{i,j} X \to X$.

Assuming the spectral sequences converge we can do the following: First, let $i_0$ be the largest integer so that $\alpha$ is the image of an element $y_0' \in \pi_m F^{i_0,-\infty} X$. Next, let $j_0$ be the largest integer so that $y_0'$ is the image of an element $y_0'' \in \pi_m F^{i_0,j_0} X$. Then the image $y_0$ of $y_0''$ in $E_1^{i_0,j_0,m-i_0-j_0} = \pi_m Gr^{i_0,j_0} X$ represents $\alpha$, and survives when we go counter-clockwise around the above square of spectral sequences.

Next let $i_1 < i_0$ denote the largest integer so that the image $y_1'$ of $y_0''$ in $\pi_m F^{i_1,j_0} X$ is the image of some element $y_1'' \in \pi_m F^{i_1,j_1} X$ with $j_1 > j_0$. We can choose $j_1$ to be maximal. Then the image $y_1$ of $y_1''$ in $\pi_m Gr^{i_1,j_1} X$ also represents $\alpha$. We can continue like this until we find a last element $y_a'' \in \pi_m F^{i_a,j_a} X$. Its image $y_a \in \pi_m Gr^{i_a,j_a} X$ represents $\alpha$ and survives when we go clockwise around the above square of spectral sequences.

The fact that the image of $y_{b-1}''$ in $\pi_m F^{i_b+1,j_{b-1}} X/F^{i_b+1,j_{b-1}+1} X$ is nonzero but maps to zero in $\pi_m F^{i_b,j_{b-1}} X/F^{i_b,j_{b-1}+1} X$ implies that there is some $z_b \in E_1^{i_b,j_{b-1},m-i_b-j_{b-1}+1} = \pi_{m+1} Gr^{i_b,j_{b-1}}$ with $d^h_{i_{b-1}-i_b}(z_b)=y_{b-1}$.

Similarly, there is some $z_b' \in E_1^{i_b,j_{b-1},m-i_b-j_{b-1}}$ with $d^v_{j_b-j_{b-1}}(z_b')=y_b$. It follows from a diagram chase that we can choose $z_b' = -z_b$. (The minus sign is not important for our purposes.) We illustrate with the following zig-zag in the case $a=2$:
\[ \xymatrix@R-18pt{
 & y_0 \\
 z_1 \ar@{|->}[ur]^-{d^h} \ar@{|->}[dr]^-{d^v} & \\
 & y_1 \\
 z_2 \ar@{|->}[ur]^-{d^h} \ar@{|->}[dr]^-{d^v} & \\
 & y_2
} \]
At one end of the zig-zag we find the representative for $\alpha$ that survives going counter-clockwise around the square of spectral sequences, at the other end we find the representative that survives going clockwise.

\subsection{The commuting square for the Tate spectrum} \label{ss:commutingsquareforTate}
Given a filtered ring $A$, we now have two compatible filtrations on $THH(A)^{tS^1}$, so by a minor reindexing of the setup in the previous section we get a commuting square of spectral sequences as follows:
\[ \xymatrix{
 \wh{E}_1^{i,j,*} = \pi_{i+j} THH(Gr A; j) \otimes P(t^{\pm 1}) \ar@{=>}[r] \ar@{=>}[d] & (\wh{E}_1')^{i+*,j} = \pi_{i+j+*} THH(Gr A;j)^{tS^1}\ar@{=>}[d] \\
 (\wh{E}_2'')^{i+j,*} = \pi_{i+j} THH(A) \otimes P(t^{\pm 1}) \ar@{=>}[r] & \pi_{i+j+*} THH(A)^{tS^1}
} \]
We have a similar commuting square of spectral sequences with relative groups everywhere computing $\pi_* THH(A,I)^{tS^1}$.

Now consider the above square for $(A,I)=(W_n(k), (p))$. Then
\[
 \wh{E}_1^{*,*,*} = THH_*(k[x]/(x^n), (x)) \otimes P(t,t^{-1}),
\]
and we understand three of the four spectral sequences completely.

Recall that we defined $z_s = t^{\lfloor s/n \rfloor} \gamma_{\lfloor s/n \rfloor}(x_n) x^{\{s/n\}}$, and that the top horizontal spectral sequence has differentials $z_s \mapsto t^{\nu_p(s)+1} \mu_0^{\nu_p(s)} z_{s-1} \sigma x$ for $n \nmid s$ and $z_s \mapsto t^{\nu_p(n)+1} \mu_0^{\nu_p(n)} z_{s-1} \sigma x$ for $n \mid s$. If we start with
\[
 \wh{E}_0^{*,*,*} = \bigoplus_{s \geq 1} k\{z_s, z_{s-1} \sigma x\} \otimes P(t,t^{-1})
\]
and interpret some of the differentials as $d_0$ differentials then this makes sense even when $p \nmid n$ and $n \mid s$.

The right hand side vertical spectral sequence collapses, because the $E_1$ term is concentrated in odd total degree.

Finally, the left hand side vertical spectral sequence has differentials which we think of as generated multiplicatively by $d_1(\mu_0)=\sigma x$, where we use the Leibniz rule $d_{r+1}(y^p)=y^{p-1} xd_r(y)$. This is only morally true, because $THH(W_n(k), (p))$ is only a non-unital ring spectrum and there is no class $\mu_0$ in $THH_2(W_n(k), (p))$. If we were considering the non-relative version of the left hand side vertical spectral sequence then this would be accurate. In any case, by Examples \ref{ex:HHofWnk} and \ref{ex:THHofW_nkrel} we know all the differentials in the left hand side vertical spectral sequence as well.

\begin{example}
For a typical example of how such a zig-zag works, let us consider the case $p=3$ and $n=4$, and let us start with the class $x_4 x \sigma x$ in total degree $3$, which survives going clockwise around the square of spectral sequences. First we use the zig-zag
\[ \xymatrix@R-18pt{
 & x_4 x \sigma x \\
 \mu_0 x_4 x \ar@{|->}[ur]^-{d^v} \ar@{|->}[dr]^-{d^h} & \\
 & t \mu_0 x_4 \sigma x \\
} \]
to replace $x_4 x \sigma x$ by the representative $t\mu_0 x_4 \sigma x$.

Next we use the zig-zag
\[ \xymatrix@R-18pt{
 & t \mu_0 x_4 \sigma x \\
 t \mu_0^2 x_4 \ar@{|->}[ur]^-{d^v} \ar@{|->}[dr]^-{d^h} & \\
 & t^3 \mu_0^4 x^2 \sigma x\sigma x \\
} \]
to replace it by the representative $t^3 \mu_0^4 x^2 \sigma x$. Here $d^h$ can be calculated using yet another zig-zag.

Finally we use the zig-zag
\[ \xymatrix@R-18pt{
 & t^3 \mu_0^4 x^2 \sigma x \\
 t^3 \mu_0^5 x^2 \ar@{|->}[ur]^-{d^v} \ar@{|->}[dr]^-{d^h} & \\
 & t^4 \mu_0^5 x \sigma x \\
} \]
to find the representative $t^4 \mu_0^4 x \sigma x$ that survives going counter-clockwise around the square of spectral sequences. (If we were considering $THH(W_4(k))$ rather than $THH(W_4(k), (3))$ we would have had a differential $d_1^v(t^4 \mu_0^6) = t^4 \mu_0^5 x \sigma x$ and our class would represent zero is $\pi_3 THH(W_4(k))^{tS^1}$.)
\end{example}

\begin{remark}
It is possible to describe the bottom horizontal spectral sequence as well. We omit the details as we will not need them.
\end{remark}

\section{A common spectral sequence for homotopy orbits and topological cyclic homology}
Suppose as usual that $A$ is a complete filtered ring or symmetric ring spectrum and let $I=F^1 A$. In addition to the usual homotopy orbit spectral sequence, we get a spectral sequence converging to $\pi_* \Sigma THH(A,I)_{hS^1}$ coming from $\Sigma THH(A,I)_{hS^1}$ being the homotopy fiber of $R : \TF(A,I) \to \TF(A,I)$. The amazing thing is that we have another spectral sequence with the same $E_1$ term converging to $\TC_*(A,I)$. The starting point is the following result.

\begin{thm}
Suppose $A$ is a complete filetered ring or symmetric ring spectrum. Then $R : \TR^{m+1}(A) \to \TR^m(A)$ sends $F^s \TR^{m+1}(A)$ to $F^{\lceil s/p \rceil} \TR^m(A)$ and $R : \TF(A) \to \TF(A)$ sends $F^s \TF(A) \to F^{\lceil s/p \rceil} \TF(A)$.
\end{thm}

\begin{proof}
We prove the case $m=1$, the general case is similar. We use the $p$-fold edgewise subdivision model of $THH$ considered in the proof of Theorem \ref{t:TRSS} above. Fixed points by the action of $C_p$ are taken spacewise, and a fixed point of a term in the colimit defining $THH^{[p]}(A; V)_q$ looks like
\[
 (a_0 \sma \ldots \sma a_q)^{\sma p} \sma v
\]
where $v \in (S^V)^{C_p}$. Now, if $a_i$ is homogeneous of filtration $|a_i|$, this is in filtration degree $p(|a_0|+\ldots+|a_q|)$. Applying $R$ replaces this by $(a_0 \sma \ldots \sma a_q) \sma v$, which has filtration degree $|a_0|+\ldots+|a_q|$.
\end{proof}

\begin{defn}
Suppose $A$ is a complete filtered ring or symmetric ring spectrum. Let $F^s \Sigma THH(A)_{hS^1}$ denote the homotopy fiber
\[
 F^s \Sigma THH(A)_{hS^1} \to F^s \TF(A) \xto{R} F^{\lceil s/p \rceil} \TF(A)
\]
and let $F^s \TC(A)$ denote the homotopy equalizer
\[
 F^s \TC(A) \to F^s \TF(A) \underset{I}{\overset{R}\rightrightarrows} F^{\lceil s/p \rceil} \TF(A).
\]
\end{defn}

\begin{thm} \label{t:TCAndhorbitSSs}
Suppose $A$ is a complete filtered ring or symmetric ring spectrum. Then there is a spectral sequence with
\begin{multline*}
 E_1^{s,t} = \ker \big( \TF_{s+t}(Gr A; s) \overset{R}{\to} \TF_{s+t}(Gr A; s/p) \big) \\
 \oplus \coker \big( \TF_{s+t+1}(Gr A; s) \overset{R}{\to} \TF_{s+t+1}(Gr A; s/p) \big)
\end{multline*}
for $s \geq 1$ and $E_1^{0,t} = 0$, converging to $\pi_* \Sigma THH(A,I)_{hS^1}$. There is another spectral sequence with the same $E_1$ term converging to $\TC_{s+t}(A,I)$.

Moreover, all differentials that multiply the filtration by a factor of less than $p$ are isomorphic and all extensions that multiply the filtration by a factor of less than $p$ are isomorphic.
\end{thm}

\begin{proof}
It is clear that there is a spectral sequence associated to the filtration of $\Sigma THH(A,I)_{hS^1}$ and another spectral sequence associated to the filtration of $\TC(A,I)$, and we can compute the $E_1$ term of the first spectral sequence using the diagram
\[ \xymatrix{
 F^{s+1} \Sigma THH(A,I)_{hS^1} \ar[r] \ar[d] & F^{s+1} \TF(A) \ar[r]^-R \ar[d] & F^{\lceil (s+1)/p \rceil} \TF(A) \ar[d] \\
 F^s \Sigma THH(A,I)_{hS^1} \ar[r] \ar[d] & F^s \TF(A) \ar[r]^-R \ar[d] & F^{\lceil s/p \rceil} \TF(A) \ar[d] \\
 Gr^s \Sigma THH(A,I)_{hS^1} \ar[r] & Gr^s \TF(A) \ar[r]^-R & Gr^{s/p} \TF(A)
} \]
A similar diagram, with $R$ replaced by $R-I$ on the top two rows, calculated the $E_1$ term of the second spectral sequence.

The last part, comparing short differentials and short extensions in the two spectral sequences, follows by Lemma \ref{l:comparepiecesofTCandhorbits} below.
\end{proof}

\begin{lemma}[Brun {\cite[Lemma 5.3]{Br01}}] \label{l:comparepiecesofTCandhorbits}
Suppose $s < t \leq ps$. Then
\[
 F^s \TC(A)/F^t \TC(A) \simeq F^s \Sigma THH(A)_{hS^1}/F^t \Sigma THH(A)_{hS^1}.
\]
\end{lemma}

This is especially useful because we can compute $\pi_* THH(W_n(k))_{hS^1}$ through a range of degrees (compare \cite[Proposition 6.4 and 7.2]{Br01}).

\begin{prop} \label{p:hoorbcalc}
For $2i \leq 2p-2$ we have
\[
 \pi_{2i} THH(W_n(k))_{hS^1} \cong W_{n(i+1)}(k)
\]
and for $2i-1 \leq 2p-3$ we have
\[
 \pi_{2i-1} THH((W_n(k))_{hS^1} = 0.
\]
\end{prop}

\begin{proof}
Recall that the homotopy orbit spectral sequence looks like
\[
 \pi_* THH(A)[t^{-1}] \Longrightarrow \pi_* THH(A)_{hS^1},
\]
and recall that through degree $2p-2$ we have $\pi_{2i} THH(W_n(k)) \cong W_n(k)$ and $\pi_{2i-1} THH(W_n(k))=0$. Hence the homotopy orbit spectral sequence collapses at the $E_2$-term through the range of degrees we consider. This then shows that $\pi_* THH(W_n(k))_{hS^1}$ has the required length over $k$.

To show that the extensions are maximally nontrivial, we consider the corresponding homotopy orbit spectral sequence with mod $p$ coefficients:
\[
 V(0)_* THH(W_n(k))[t^{-1}] \Longrightarrow V(0)_* THH(W_n(k))_{hS^1}.
\]
Let $\beta_n$ denote the element in $V(0)_1 THH(W_n(k))$ which is coming from $p^{n-1} \in THH_0(W_n(k)) \cong W_n(k)$. Then we have an immediate differential
\[
 d_2(t^{-1}) \doteq \beta_n
\]
and it follows that we have a differential
\[
 d_2(t^{-i}) \doteq t^{-i+1} \beta_n
\]
for all $i \leq p-1$. This implies that $V(0)_{2i} THH(W_n(k))_{hS^1} \cong k$ for $2i \leq 2p-2$, and it follows that the extensions are maximally nontrivial.

The result now follows because the only maximally nontrivial extension of $W_{ni}(k)$ by $W_n(k)$ is $W_{n(i+1)}(k)$. (Recall that all the identifications with Witt vectors are additive only; no multiplicative structure is implied.)
\end{proof}

\begin{cor}
For $2i-1 \leq 2p-3$ we have
\[
 \pi_{2i-1} F^1 \Sigma THH(W_n(k))_{hS^1} \cong W_{(n-1)i}(k)
\]
and for $2i \leq 2p-2$ we have
\[
 \pi_{2i} F^1 \Sigma THH(W_n(k))_{hS^1} = 0.
\]
\end{cor}

\begin{proof}[Proof of Theorem \ref{t:mainorder}]
Because $\TF_*(W_n(k), (p))$ is concentrated in odd total degree it follows that in this case all the differentials in the above spectral sequence converging to $\TC_*(W_n(k), (p))$ go from odd to even total degree. The result now follows by a counting argument.
\end{proof}

If instead we use the non-relative $K$-theory spectrum $K(W_n(k))$, we pick up an extra $K(k)$ and we get the following.

\begin{cor} \label{c:ordernonrel}
Suppose $k=\bF_q$ is a finite field with $q$ elements. Then
\[
 \frac{|K_{2i-1}(W_n(k))|}{|K_{2i-2}(W_n(k))|} = q^{(n-1)i}(q^i-1)
\]
for all $i \geq 2$.
\end{cor}

Finally we can prove the last main result.

\begin{proof}[Proof of Theorem \ref{t:mainlowdeg}]
We compare the differentials and extensions in the spectral sequence converging to $\pi_* \Sigma THH(W_n(k), (p))_{hS^1}$ to the differentials and extensions in the spectral sequence converging to $\TC_*(W_n(k), (p))$, using Theorem \ref{t:TCAndhorbitSSs}.

To compute differentials we need to calculate $R$. To do this, we use the commutative diagram
\[ \xymatrix{
 \TF_*(W_n(k), (p)) \ar[r]^-R \ar[d]^\Gamma & \TF_*(W_n(k), (p)) \ar[d]^{\wh{\Gamma}} \\
 \pi_* THH(W_n(k), (p))^{hS^1} \ar[r]^-{R^h} & \pi_* THH(W_n(k), (p))^{tS^1}
} \]
and the observation (Theorems \ref{t:Gamma_injective} and \ref{t:Gammahat_injective}) that $\Gamma$ and $\wh{\Gamma}$ are injective. Given an element in $\TF_*(W_n(k), (p))$, we need to be able to find a representative for its image under $\wh{\Gamma}$ in the Tate spectral sequence; for this we use the ``commutative square''
\[ \xymatrix{
 \wh{E}_1^{*,*,*} = THH_*(k[x]/(x^n), (x)) \otimes P(t,t^{-1}) \ar@{=>}[r] \ar@{=>}[d] & \pi_* THH(k[x]/(x^n), (x))^{tS^1} \ar@{=>}[d] \\
 THH_*(W_n(k), (p)) \otimes P(t,t^{-1}) \ar@{=>}[r] & \pi_* THH(W_n(k), (p))^{tS^1}
} \]
of spectral sequences discussed in Section \ref{s:commutingsquareofSSs} above.

We use the collapsing spectral sequence
\[
 E_1^{*,*} = \TF_*(k[x]/(x^n), (x))^{hS^1} \Longrightarrow \TF_*(W_n(k), (p))
\]
and the injective map $\Gamma : \TF_*(k[x]/(x^n), (x)) \to \pi_* THH(k[x]/(x^n), (x))^{hS^1}$ to name elements in $\TF_*(W_n(k), (p))$.

Recall that we defined
\[
 y_s = \mu_0^{-\lfloor s/n \rfloor} \gamma_{\lfloor s/n \rfloor}(x_n) x^{\{s/n\}}.
\]
We are interested in the class $\mu_0^a y_{s-1} \sigma x$, which lives in $\TF_{2a+1}(k[x]/(x^n); s)$ for $a \geq \lfloor \frac{s-1}{n} \rfloor$. (If $p \nmid n$ and $n \mid s$ then this class is zero.) Through the range of degrees we are interested in we have $a \leq p-1$.

Also recall that we defined
\[
 z_s = t^{\lfloor s/n \rfloor} \gamma_{\lfloor s/n \rfloor}(x_n) x^{\{s/n\}}.
\]
Then $t^b z_{s-1} \sigma x$ lives in the $E_2$ term of the Tate spectral sequence converging to $\pi_* THH(k[x]/(x^n); s)^{tS^1}$ for all $b$.

The heart of the argument is the following. Even though the class $t^b z_{s-1} \sigma x$, interpreted as lying in $\wh{E}_1^{*,*,*}$, is killed by a differential in the vertical spectral sequence converging to $THH_*(W_n(k), (p)) \otimes P(t,t^{-1})$, we can use the above commutative square of spectral sequences to find a different representative for this class in $\wh{E}_1^{*,*,*}$.

Start with a non-zero class $\mu_0^a y_{s-1} \sigma x \in \TF_{2a+1}(k[x]/(x^n);s)$ with $\lfloor \frac{s-1}{n} \rfloor \leq a \leq p-1$. Then
\[
 \wh{\Gamma}(\mu_0^a y_{s-1} \sigma x) = t^{-a} z_{ps-1} \sigma x.
\]
If the power of $t$ in
\[
 t^{-a} z_{ps-1} \sigma x = t^{\lfloor (ps-1)/n \rfloor -a} \gamma_{\lfloor (ps-1)/n \rfloor}(x_n) x^{\{(ps-1)/n\}} \sigma x
\]
is non-negative then there is nothing to do: $t^{-a} z_{ps-1} \sigma x$ lies in the restriction of the Tate spectral sequence for $(W_n(k), (p))$ to the second quadrant and hence $t^{-a} z_{ps-1} \sigma x$ is in the image of the map $R^h$. We conclude that $\mu_0^a y_{s-1} \sigma x$ is $R$ of a class in filtration $ps$. (We already knew that, because in this case $R : \TF_{2a+1}(k[x]/(x^n); ps) \to \TF_{2a+1}(k[x]/(x^n); s)$ is an isomorphism.)

If the power of $t$ in $t^{-a} z_{ps-1} \sigma x$ is negative we can use the zig-zag
\[ \xymatrix@R-18pt{
 & t^{-a} z_{ps-1} \sigma x \\
t^{-a} \mu_0 z_{ps-1} \ar@{|->}[ur]^-{d^v} \ar@{|->}[dr]^-{d^h} & \\
 & t^{-a+1} \mu_0 z_{ps-2} \sigma x
} \]
to replace $t^{-a} z_{ps-1} \sigma x$ with the representative $t^{-a+1} \mu_0 z_{ps-2} \sigma x$. If $p \nmid n$ and $n \mid ps-1$ then the horizontal differential on $t^{-1} \mu_0 z_{ps-1}$ is not given by this formula but rather by applying another zig-zag
\[ \xymatrix@R-18pt{
 & t^{-a+1} \mu_0 z_{ps-2} \sigma x \\
t^{-a+1} \mu_0^2 z_{ps-2} \ar@{|->}[ur]^-{d^v} \ar@{|->}[dr]^-{d^h} & \\
 & t^{-a+2} \mu_0^2 z_{ps-3} \sigma x
} \]
but the end result is the same.

We can continue like this until we obtain a representative $t^{-a+b} \mu_0^b z_{ps-b-1} \sigma x$ with a non-negative power of $t$. By our assumptions this is always possible, and it happens for some $b \leq p-1$. (Through this range of degrees we never encounter $\mu_0^b$ for $p \mid b$, so the differential on $\mu_0^b$ is always as simple as possible. We also never encounter $z_c$ for $p \mid c$, so the differential on $z_c$ is also as simple as possible.) We conclude that this class is in the image of $R^h$, and that $\mu_0^a y_{s-1} \sigma x$ is $R$ of a class in filtration $ps-b$.

A similar discussion shows that multiples of $\mu_0^a y_{s-1} \sigma x$ are in the image of $R$, starting with the fact that $\wh{\Gamma}(p^{a'} \mu_0^a y_{s-1} \sigma x)$ is represented by $t^{-a+a'} \mu_0^{a'} z_{ps-1} \sigma x$ in the Tate spectral sequence converging to $\pi_* THH(k[x]/x^n, (x))^{tS^1}$.

From the above discussion we see that through degree $2p-1$ the only ``long'' differential in the spectral sequence converging to $\pi_* \Sigma THH(W_n(k), (p))_{hS^1}$ is the one coming from $R(\mu_0^{p-1} \sigma x) = \mu_0^{p-1} \sigma x$. With our conventions, $R$ is given by the inverse Frobenius (see Section \ref{ss:TCk}) and is an isomorphism
\[
 k \cong E_1^{1,2p-2} = E_p^{1,2p-2} \xto{d_p=R} E_p^{p,p-2} = E_1^{p,p-2} \cong k.
\]

In the spectral sequence for topological cyclic homology we have
\[
 k \cong E_1^{1,2p-2} = E_p^{1,2p-2} \xto{d_p=R-1} E_p^{p,p-2} = E_1^{p,p-2} \cong k.
\]
The kernel of $d_p$ in this particular bidegree is $\bZ/p$ and the cokernel is $\coker(\phi - 1) : k \to k$.

This proves Theorem \ref{t:mainlowdeg} up to extensions. It follows from Theorem \ref{t:TCAndhorbitSSs} and Proposition \ref{p:hoorbcalc} that all the extensions except one are maximally nontrivial. The one extension we cannot calculate in this way goes from filtration $1$ to filtration $3$ (or, in some cases, to filtration $4$ or $5$) in degree $2p-3$, because such an extension only makes sense after the $d_{p-1}$ differential happens, and $F^1 \TC(W_n(k))/F^p \TC(W_n(k))$ is just outside the range where Lemma \ref{l:comparepiecesofTCandhorbits} applies. The picture below shows what the spectral sequence looks like in this degree in the case $p=5$, $n=5$; the dashed line is the extensions we cannot calculate using Theorem \ref{t:TCAndhorbitSSs}.

\begin{center}
 \includegraphics[scale=0.4]{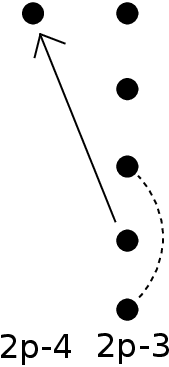}
\end{center}

The fact that this particular extensions is as claimed follows by comparing with $\TC_{2p-3}(W(k), (p))$.
\end{proof}

% \bibliographystyle{plain}
% \bibliography{b}

\end{document}